\input amstex
\documentstyle{amsppt}
\vsize=8in \hsize 6.6 in 
\loadbold \topmatter
\title Flag Foliations Functionals. \\ The Hopf Hypothesi. \endtitle
\author Valery Marenich \endauthor

\address Kalmar, Sweden\endaddress \email valery.marenich\@ gmail.com \endemail
 \subjclass 58, 53C \endsubjclass
\date{August 24 of 2009}\enddate
\abstract We study the Riemannian manifold $(M,g)$ diffeomorphic to the direct product $\bar M^4=S^2\times S^2$
of two two-spheres  by some diffeomorphism $f: S^2\times S^2 \to M$. The flag foliation $(\bar\Upsilon, \bar
\Cal F)$ on  $\bar M^4$ by the push-forward construction by $f$ provides the flag foliation $(\Upsilon, \Cal F)$
on $M^4$. We define some functional ${\Cal H}$ on leaves of the foliation ${\Cal F}$. Applying the Nash process
we construct the smooth $f^*$ such that ${\Cal H}$ vanishes on some, so called, "silent" leaf ${\Cal F^*}$,
implying the affirmative solution to the Hopf conjecture on the non-existence of the Riemannian metric $g$ on
$S^2\times S^2$ with positive sectional curvature.
\endabstract
\endtopmatter

\document

\head Introduction \endhead

The direct product  $\bar M^4=S^2\times S^2$ contains many totally geodesic tori. We may choose one of them
${\bar\Cal F^*}$ and construct the torus action on $\bar M^4$ such the set of its orbits will give the
(singular) foliation ${\bar\Cal F}_{\alpha,\beta}$ of $\bar M^4$ with ${\bar\Cal F}^*$ being one of the leaves.
Starting to change the metric on $\bar M^4$ we then observe what happens with this foliation and ${\bar\Cal
F}^*$. We can not expect that the totally geodesic torus "survives" this change, but may hope that - possibly
after some correction by a diffeomorphism of $M^4$ - still there will be in arbitrary $(M^4,g)$ diffeomorphic to
$\bar M^4$ some "special" torus ${\Cal F}^*$ flat enough to guarantee the existence of zero curvature somewhere.
For instance, due to the Gauss equation and the Gauss-Bonnet theorem this happens if in each point of ${\Cal
F}^*$ there exists a direction in which all second forms of ${\Cal F}^*$ vanish. Assume that these are
directions tangent to yet another one-dimensional foliation $\Upsilon$ on $M^4$ defined in a similar way as
$\Cal F$: take the orbits of the diagonal $S^1$-action of the $S^1\times S^1$-torus action on $S^2\times S^2$.
The "flatness" of the each particular torus ${\Cal F}_{\alpha,\beta}$  then will be measured by the maximum
rotation of its tangent planes along our wiring $\Upsilon$, see the definition of the functional ${\Cal H}({\Cal
F}_{\alpha,\beta})$ below. If ${\Cal H}({\Cal F}_{\alpha,\beta})=0$ - we say that ${\Cal F}_{\alpha,\beta}$ is
$\Upsilon$-parallel\footnote{or $\Upsilon$-holomorphic, see Appendix for explanations.} - then the tangent
planes are parallel along all curves $\gamma\in\Upsilon$ and the Hopf conjecture follows. Connecting an
arbitrary $(M^4,g)$ with the standard direct product $(\bar M^4,\bar g)$ by a family $(M^4,g(\tau)),
0\leq\tau\leq 1$, where $g(\tau)$ is a smooth family of Riemannian metrics, we know that for $\tau=0$ such torus
exists. How we can be sure that it does not disappear? We associate to every flag foliation $(\Upsilon, \Cal F)$
a vector field ${\Cal W}(\alpha,\beta)$ on the moduli space $X^2$ of the foliation ${\Cal F}$ on $M^4$ (in our
case this is the square $[-\pi/2,\pi/2]\times[-\pi/2,\pi/2]$) such that $\Upsilon$-parallel leaves would be
zeroes of ${\Cal W}$. Then we check that $\Cal W$ has one zero of the index one when the metric on $M^4$ of the
standard direct product; and that during the deformation $g(\tau)$ the vector field ${\Cal W}$ never vanish on
the boundary of $X^2$ (non-collapsing). This provides us for every $\tau$ with some torus ${\Cal F}_{\alpha,\beta}$
such that ${\Cal W}({\alpha,\beta})=0$, and the last condition makes the Nash minimizing process possible, resulting
in some $\Upsilon$-parallel ${\Cal F^*}$.

Remark, that the given metric $g$ is the limit of $g(\tau)$ when $\tau\to 1$, where the family of metrics is
chosen to satisfy some "standard regularity assumptions". Therefore, we do not claim that the "silent" torus
exists in the arbitrary given $(M^4,g)$. Instead the family of planes of zero curvature tangent to the silent
tori in $(M^4,g(\tau))$ will have some limit, due to compactness of $M^4$ - some two-dimensional direction of
zero-curvature in $(M^4,g)$, which infer the Hopf conjecture.

\medskip

\head 1. Flag Foliations \endhead

Take two two-spheres $\bar S^2_+$ and $\bar S^2_-$ in a three-dimensional  Euclidean space, and denote by $\bar
l_\alpha(\phi)$ and $\bar l_\beta(\psi)$  two families of "parallels" on $\bar S^2_+$ and $\bar S^2_-$ which are
correspondingly: orbits of an $S^1$-action by rotation on the angle $\phi$ about the axe in space going through
the north and south poles $N_+$ and $S_+$ of the sphere $\bar S^2_+$, and similar $S^1$-action by rotation on
the angle $\psi$ about the axe in space going through the north and south poles $N_-$ and $S_-$ of the sphere
$\bar S^2_-$. Parameters $\alpha$ and $\beta$ are "altitudes" and vary from $-\pi/2$ (south poles) to $\pi/2$
(north poles). On the direct product $\bar M^4=\bar S^2_+\times \bar S^2_-$ the torus action $S^1(\phi)\times
S^1(\psi)$ is defined with orbits ${\bar{\Cal F}}^2_{(\alpha,\beta)}(\phi,\psi)$, where ${\bar{\Cal
F}}^2_{(\alpha,\beta)}=\bar l_\alpha(\phi)\times \bar l_\beta(\psi)$, and $(\alpha,\beta,\phi,\psi)$ may be
considered as "coordinates" of the points in $\bar M^4$. For $-\pi/2<\alpha,\beta<\pi/2$ these orbits are tori
which degenerates to a circle when one of $\alpha$ or $\beta$ converge to $\pm \pi/2$, or to one of the four
fixed points of the torus action when both $\alpha,\beta$ converge by an absolute value to $\pi/2$. The square
$X^2=[-\pi/2,\pi/2]\times [-\pi/2,\pi/2]$ of all pairs $(\alpha,\beta)$ is the space of all orbits of the torus
action, which we call the moduli space. This family ${\bar{\Cal F}}^2_{(\alpha,\beta)}$ is a singular foliation
of $\bar M^4$ which we denote by $\bar\Cal F$. Taking $\phi=\psi=t$ we obtain the $S^1$-action on $\bar M^4$
with orbits=circles which provide wirings of every tori ${\bar{\Cal F}}^2_{\alpha,\beta}$: if
$\bar\gamma_{(\alpha,\beta,\phi,\psi)}(t)$ such an orbit issuing from the point with "coordinates"
$(\alpha,\beta,\phi,\psi)$ then all the points $(\alpha,\beta,\phi+t,\psi+t)$ belong to the same orbit of this
$S^1$-action, and the circle $\bar \gamma_{(\alpha,\beta,\phi,\psi)}$ has the homotopy class $(1,1)$ in the
fundamental group $\pi_1({\bar{\Cal F}}^2_{(\alpha,\beta)})$ of the torus ${\bar{\Cal
F}}^2_{(\alpha,\beta)}=\bar l_\alpha(\phi)\times \bar l_\beta(\psi)$ generated by $\bar l_\alpha$ and $\bar
l_\beta$ ("parallels" and "meridians"). The family of these circles $\bar \gamma_{(\alpha,\beta,\phi,\psi)}(t)$
we denote by $\bar\Upsilon$, it is a singular one-dimensional foliation of $\bar M^4$. At every "regular" point
$\bar P$ of $\bar M$ - except with $|\alpha|=|\beta|=\pi/2$ - we have two vector spaces in the tangent space
$T_{\bar P}\bar M$: two-dimensional tangent to the foliation $\bar\Cal F$ and one-dimensional subspace tangent
to the foliation $\bar\Upsilon$. These vector spaces give a flag in $T_{\bar P}\bar M$, and we call the pair
$(\bar\Upsilon, \bar \Cal F)$ - {\bf the flag foliation}.

\medskip

Next: an arbitrary diffeomorphism $f:\bar M^4\to M^4$ provides $M^4$ by a  push-forward construction with the
pair of (singular) foliations which we denote by $(\Upsilon,\Cal F)$ and call {\bf the flag
foliation}.\footnote{It appears, that flag foliations are, in a sense, dual to open book decompositions, see
Appendix~2.} Leaves of these foliations going through the point $P$ of $M^4$ we denote by $\gamma_P(t)$ and
${\Cal F}_P$ correspondingly. Another notations are:
$$
{\Cal F}^2_{(\alpha,\beta)}=f(\bar{\Cal F}^2_{(\alpha,\beta)})\qquad \text{ and } \qquad \gamma_{(\alpha,\beta,\psi,\phi)}(t)=f(\bar\gamma_{(\alpha,\beta,\phi,\psi)}(t)).
$$
By $I_\gamma(t)$ we denote the operator of a parallel translation in $M^4$  along $\gamma_P$ from the point
$\gamma_P(t)$ to the point $P$, assuming that parameter $t$ on $\gamma_P$ is such that $P=\gamma_P(0)$. To the
tangent bundle $TM^4$ we associate two other bundles - $SO$ - the bundle with the fiber $SO(4)$ over $P$ of
orientation preserving isometries of $T_PM^4$, and the (Grassmann) bundle $G$ of two-planes with the fiber
$SO(4)/S^1\times S^1$. Finely, let $T(t)$ be the tangent plane to $F_P$ at the point $\gamma_P(t)$, and ${\Cal
H}(\gamma_P)$ denotes the length of the curve $L(t)=I_\gamma(t) (T(t))$ of two-planes in the Grassmann fiber
$G_P$ over $P$ obtained by parallel translations of $T(t)$ along $\gamma_P$ to $P$. Then
$$
{\Cal H}({\Cal F}_{\alpha,\beta})=\max\limits_{\gamma_P\subset {\Cal F}_{\alpha,\beta}} {\Cal H}(\gamma)
\qquad \text{ and } \qquad {\Cal H}(M^4,g,f)=\min\limits_{{\Cal F}_{\alpha,\beta}} {\Cal H}({\Cal F}_{\alpha,\beta}). \tag 1.1
$$

Considering one torus ${\Cal F}_{\alpha^*,\beta^*}$ we usually employ the  pair $\{\psi,t\}$ as coordinates on
it, and denote $\gamma_{(\alpha,\beta,\phi,\psi)}(t)$ simply by $\gamma_\psi(t)$. Then $T_\psi(t)$ and
$L_\psi(t)$ denote correspondingly the family of tangent planes to the torus along $\gamma_\psi(t)$ and the
parallel transport of this family by $I_{\gamma_\psi}(t)$ to the initial point $P_\psi=\gamma_\psi(0)$.

\medskip

Besides the functional ${\Cal H}({\Cal F}_{\alpha,\beta})$ on leaves of the  foliation ${\Cal F}$ the second
most important object we consider is the following "vector" field ${\Cal W}$\footnote{a field of directions
actually, since its coordinates $\Phi,\Psi$ are defined only up to the simultaneous change of a sign} on the
moduli space $X$ of tori leaves ${\Cal F}_{\alpha,\beta}$, which is defined as follows. Each parallel
translation $I_\gamma$ acts on its invariant planes $\Pi$ and $N$ by rotations on angles
$\Phi(\alpha,\beta,\psi)$ and $\Psi(\alpha,\beta,\psi)$. Then
$$
{\Cal W}(\alpha,\beta, \psi)=(\Phi(\alpha,\beta,\psi),\Psi(\alpha,\beta,\psi))
$$
is "the profile curve" of ${\Cal F}_{\alpha,\beta}$ and
$$
{\Cal W}(\alpha,\beta)=\int\limits_0^{2\pi}(\Phi(\alpha,\beta,\psi),\Psi(\alpha,\beta,\psi)) d\psi. \tag 1.2
$$
We may assume that this vector field is "good", i.e.; it has only  non-degenerated zeroes inside the moduli
space $X^2$. Denote by $ind(x)$ the index of zero $x$ of ${\Cal W}$ and set $ind(x)=0$ for non-zeros. Finely,
$$
{\Cal He}(M,g,f)= \sum\limits_{x\in X^2} {ind(x(\alpha,\beta))} {\Cal H}({\Cal F}_{\alpha,\beta}), \tag 1.3
$$
where $x(\alpha,\beta)$ represent ${\Cal F}_{\alpha,\beta}$ in the  moduli space $X^2$.\footnote{Note, that
(1.3) looks like the Euler characteristic of some homology group related to $M^4$; i.e., the Euler co-cycle on a
fundamental class of $X^2$ weighted by ${\Cal H}$.}

Below we show that for $M^4$ diffeomorphic to $S^2\times S^2$ the  index of the vector field ${\Cal W}$ along
the boundary of $X$ equals one, or it always has at least one zero inside $X^2$. In other cases, say for $S^4$ -
the suspension over $S^3$ with a foliation by the Clifford tori; the field ${\Cal W}$ has no zeroes. Note that,
slightly changing if necessary,\footnote{or introducing small "fake" holonomies} we may assume that under a
given smooth family $f_\delta$ of diffeomorphisms zeros of ${\Cal W}$ behave in a general way: for all but a
finite number $\delta_i$ of instances of $\delta$ they are isolated and non-degenerated; while at each of the
moment $\delta_i$ the pair of zeros of opposite indexes $\pm 1$ is born or annihilated. This implies that our
functional ${\Cal He}(M,g,f)$ is continuous and - with a possible exception of the set of $\delta_i$ - is
differentiable on $\delta$.

Our Nash process is diminishing ${\Cal H}({\Cal F}_{\alpha,\beta})$ on those leaves ${\Cal F}_{\alpha,\beta}$
which contribution to (1.2) is +1. We make sure that ${\Cal H}({\Cal F}_{\alpha^-,\beta^-})$ on leaves with
contribution -1, which may be born together with a companion ${\Cal F}_{\alpha^+,\beta^+}$ of index +1, is not
less than ${\Cal H}({\Cal F}_{\alpha^+,\beta^+})$ for this companion - which is possible since they are born
with equal ${\Cal H}$'s. Since ${\Cal H}$ of the two annihilating zeroes of ${\Cal W}$ of opposite indexes are
equal too this ensures the surviving of the zero ${\Cal F}_{\alpha^*,\beta^*}$ of ${\Cal W}$ with index +1 and
minimal $\Cal H$ for all $\tau$ and $\delta$ (since the zero of index +1 may annihilate only with - popping up
on ${\Cal H}$-hight scale - zero of index -1, which index +1 companion has smaller $\Cal H$; and which we know
exists at the beginning of the whole story). During this process ${\Cal H}(M^4,g,f)$ may jump, but only down,
and achieves its minimum zero on $f^*$.

\medskip

All in the text below is about the holonomy operator, so the following Ambrose-Singer formula for the
deformation of the holonomy operator plays a crucial role. If $\gamma_\delta(t)$ is the smooth family of closed
curves, $W_\delta(0)$ is some vector field along $\gamma_\delta(0)$, and
$W_\delta(t)=I^{-1}_{\gamma_\delta(t)}(W_\delta(0))$ the parallel translation of $W_\delta(0)$ along
$\gamma_\delta$ from $\gamma_\delta(0)$ to $\gamma_\delta(t)$, then
$$
\nabla_{{\partial}/{\partial\delta}} (W_\delta(2\pi) - W_\delta(0)) = \int\limits_0^{2\pi} I_{\gamma_\delta}(t) (R(V_\delta(t),\dot\gamma_\delta(t))W_\delta(t), \tag 1.4
$$
where $V_\delta(t)$ is the variation field $\partial\gamma_\delta(t)/\partial\delta$.

\medskip

\head 2. Variations of ${\Cal H}(\gamma)$ \endhead

We study the behavior of our functional under the variations ("corrections": "silencing" and "tunings") of the
diffeomorphism $f:\bar M^4\to M^4$ by taking compositions from the left $f_\delta=\phi_\delta \circ f$ with
one-parameter families of diffeomorphisms $\phi_\delta$ of $M^4$.\footnote{Another equivalent way is by taking
compositions from the right $f_\delta=f\circ \bar\phi_\delta$ with one-parameter families of diffeomorphisms
$\bar\phi_\delta$ of $\bar M^4$. But then tame estimates follow after longer arguments.} Families $\phi_\delta$
will be generated by vector fields $V$ on $M^4$. Below notations with subscript $\delta$ denote corresponding
objects under the family of deformations $f_\delta$; while $\delta X$, $\delta Y$ or $\delta T(t)$ (and others
with pre-scrip $\delta$) will denote the finite change (usually of order $O(\delta)$) of the corresponding
linear object.

Let $\Upsilon_{\alpha,\beta}$ be the set of those $\gamma\subset {\Cal F}_{\alpha,\beta}$ where ${\Cal
H}(\gamma)={\Cal H}({\Cal F}_{\alpha,\beta})$, which we call "extremals", i.e., those $\gamma^*$ where our
functional attains its maximum on the leaf. For all other $\gamma$ from the same leaf it holds ${\Cal H}(\gamma)
< {\Cal H}(\gamma^*)$. The union of all $\gamma\in\Upsilon_{\alpha,\beta}$ with ${\Cal H}(\gamma)< {\Cal
H}(\gamma^*)- \epsilon$ is an open set $U_{\alpha,\beta}(\epsilon,f)$ converging to the set of extremals as
$\epsilon\to 0$. Since ${\Cal H}(\gamma)$ continuously depends on $\gamma$ and $f$, the small variations of $f$
inside $U(\epsilon,f)$ does not affect ${\Cal H}({\Cal F}_{\alpha,\beta})$, so that
$$
{{\partial}\over{\partial\delta}}{\Cal H}({\Cal F}_{\alpha,\beta}, f_\delta)_{|\delta=0} =
\max_{\Upsilon_{\alpha,\beta}} \{ {{\partial}\over{\partial\delta}} {\Cal H}(\gamma_\delta)_{|\delta=0} \}, \tag 2.1
$$
i.e., in order to diminish the functional we should construct such corrections of $f$ which diminish its value
on extremals.

The  only necessary restriction on our deformations is that they are "tame", or that they are chosen in such a
way that their superposition converges to some smooth limit.

In the next section~3  we construct the curve $f_\delta$ (the Nash process) in the space $Diff(\bar M^4,M^4)$ of
diffeomorphisms from $\bar M^4$ to $M^4$ converging to the smooth solution $f^*$ of the minimizing problem for
${\Cal H(F)}$. To prove convergence we employ the theory of tame spaces and maps developed in [H]. Here we
recall some facts of this theory we need below, construct elementary deformations and verify tame estimates for
them.

\medskip

\subhead 2.1. The Tame Structure on $Diff(\bar M^4, M^4)$ \endsubhead

It is known that the group of smooth diffeomorphisms  $Diff(M)$ of an arbitrary smooth manifold has the
structure of a  Fr\'{e}chet manifold and is a Fr\'{e}chet Lie group. Its Lie algebra, the tangent space at the
identity map is the space ${\Cal C}^\infty(M)$ of smooth vector fields on $M$ with the bracket operation - the
usual Lie bracket of vector fields. In a similar way, the set of all diffeomorphisms $Diff(\bar M^4, M^4)$ is
also the Fr\'{e}chet manifold with the topology defined by the family of $C^n$-norms. For the given
diffeomorphism $f:\bar M^4\to M^4$ the tangent space to this manifold may be identified with ${\Cal
C}^\infty(M^4)$ or with ${\Cal C}^\infty(\bar M^4)$. These are the spaces of vector fields $V$ on $M^4$ and
$\bar V$ on $\bar M^4$ correspondingly. In the construction of the Nash process below we employ the smoothing
operators $S_\theta$ on these spaces. One way to construct them is taking convolutions with smooth functions of
small support and orthogonal to all polynomials (in some fixed system of local coordinates) up to some degree;
see [G], the chapter 2.2.3 "The Nash (Newton-Moser) Process", pages 121--129. Another (equivalent) way is to use
the theory of the tame spaces and maps constructed in [H] according to which the tame structure on ${\Cal
C}^\infty(\bar M^4)$ may be introduced as follows:

Embed $\bar M^4$ into an open ball $\bar B^6$ of the Euclidean spaces $\bar M^4 \hookrightarrow \bar B^6 \subset
R^6$, and denote by ${\Cal C}_0^\infty=C_0^\infty(R^6,R^6)$ the Fr\'{e}chet space of smooth functions $\bar
W:R^6\to R^6$ which all derivatives tend to zero at infinity with the topology given by a family of $C^n$-norms:
$$
\|\bar W\|_n=\sup_{|\alpha|=n} \sup_{\bar M} |D^\alpha \bar W|.
$$
Then we  find a continuous extension operator $\varepsilon: {\Cal C}^\infty(\bar M^4) \to {\Cal C}_0^\infty$
such that
$$
\|\varepsilon(\bar V)\|_n \leq C \| \bar V\|_{C^n}.
$$
Say, choose tubular neighborhoods of $\bar M^4$ in $R^6$, extend $\bar V$ to be constant along the fibres, and
cut off with a smooth bump function. If $i$ - inclusion, $\rho$- restriction then
$$
\|i \bar V\|_n \leq C\| \bar V \|_n \qquad \text{ and } \qquad \| \rho \bar V \|_n \leq C \| \bar V \|_n
$$
and we have
$$
{\Cal C}^\infty(\bar M^4) {\overset{\varepsilon}\to{\longrightarrow}} C_0^\infty(\bar B^4){\overset{ i}\to{\longrightarrow}} {\Cal C}_0^\infty
{\overset{\rho}\to{\longrightarrow}}{\Cal C}^\infty(\bar M^4),
$$
where  the long composition above is the identity map. With the help of the Fourier transform ${\Cal Fr}$ the
map $i$ may be factored through the Fr\'{e}chet spaces $L_1^\infty (R^6, R^6, d\xi , log(1+|\xi|))$ with the
family of norms
$$
\|\bar W\|_n = \int e^{n log(1+|\xi|)} \|\bar W(\xi)\| d\xi .
$$
as
$$
C_0^\infty(\bar B^4){\overset{\Cal Fr}\to{\longrightarrow}} L_1^\infty(R^6, R^6, d\xi, log(1+|\xi|)) {\overset{ {\Cal Fr}^{-1}}\to{\longrightarrow}}
{\Cal C}_0^\infty
$$
As in [H] the maps ${\Cal Fr}$ and ${\Cal Fr}^{-1}$ are tame, i.e., satisfy the following estimates:\footnote{"the lost of differentiability is uniformly bounded"}
$$
\|Fr^{-1}(\bar W)\|_n \leq C\| \bar W\|_n \qquad \text{ and } \qquad \| Fr(\bar W)\|_n \leq C\|\bar W \|_{n+r}
$$
for some  $r$, which shows that ${\Cal C}^\infty(\bar M^4)$ with a family of $C^n$-norms as above is tame, see
again [H], pages 136-138 for details. In the tame space, which can be identified with the space $\Sigma(B)$ of
sequences $\{w_l\}$ in some Banach space $B$ with $\|\{w_l\}\|_n=\sum e^{nl}\|w_l\|_B$; the smoothing operators
have very simple form:
$$
\bar S_\theta (\{w_k\})=\{s(\theta-l) w_l\},
$$
where $s(\theta)$  is a smooth function which equals $1$ for $\theta\leq 0$, $0$ for $1\leq \theta$ and $0\leq
s(\theta)\leq 1$ in between for $0\leq \theta\leq 1$.\footnote{i.e., the smoothing is the cutting out higher
harmonics from the function.} In the same way we construct the smoothing operators $S_\theta$ on the space
${\Cal C}^\infty(M^4)$ of smooth vector fields on $M^4$.

Both smoothing operators satisfy the following estimates:
$$
\|S_\theta(V)\|_n \leq C e^{n-m}\theta \|V\|_m \qquad \text{ and } \qquad \|(I-S_\theta)(V)\|_m \leq C e^{m-n}\theta \|V\|_n, \tag 2.1.1
$$
see [H], page 176.

Similarly, using some isometric embedding $M^4 \hookrightarrow R^N$ the tame structure and smoothing operators
are introduced for $M^4$.

\subhead 2.1.2. Smoothing deformations \endsubhead

Our deformations will be localized near some leaf ${\Cal F}_{\alpha^*,\beta^*}$ - in its open
"$\epsilon$-neighborhood" ${\Cal O}_{\alpha^*,\beta^*,\epsilon}$ in coordinates $(\alpha, \beta, \psi,t,)$;
i.e.,
$$
{\Cal O}_{\alpha^*,\beta^*,\epsilon} = \{\gamma_\psi(t) | \gamma_\psi \subset {\Cal F}_{\alpha,\beta} \quad \text{ where } \quad (\alpha-\alpha^*)^2+(\beta-\beta^*)^2 < \epsilon^2\}.
$$

It is important not to forget that in order to be smoothing operators as above our corrections should be written
in corresponding coordinates in the Euclidean space of the isometric embedding $M^4 \hookrightarrow R^N$. Nevertheless, for
our convenience, we write them in our coordinates $(\alpha,\beta,\phi,\psi)$ since the necessary "translation" to
coordinates in the Euclidean space $R^N$ of the isometric embedding is obvious.

As above, $\Pi_\psi$ and $N_\psi$ denote invariant planes of the parallel transport $I_{\gamma_\psi}$ along
closed $\gamma_\psi$. Let also  $E_i(\psi), i=1,2$ and $E_j(\psi), j=3,4$ be the orthonormal bases at
$\gamma_\psi(0)$ of $\Pi_\psi$ and $N_\psi$ correspondingly. If we assume for the moment that that $E_i(\psi)$
are smooth vector fields along closed meridian $l(\psi)$ on ${\Cal F}_{\alpha^*,\beta^*}$, then the functions
$x^i$ below are smooth periodic on $\psi$. In general, $\Pi_\psi$ and $N_\psi$ are uniquely (correctly) defined
and depend smoothly on $\psi$ only for those $\psi$ where $I_{\gamma_\psi}\not= id$, i.e., on some disjoint
union of open subintervals $(\psi^-_m,\psi^+_m)$ of the circle $0\leq\psi\leq 2\pi$. The same is true for
$E_i(\psi)$. On the closed intervals $[\psi^+_m,\psi^-_{m+1}]$ the vectors fields $E_i(\psi)$ may be extended
arbitrarily - we choose them to provide the smooth if possible; i.e. when $\psi^+_m \not= \psi^-_{m+1}$. As the
result we will have $E_i(\psi,t)$ and $x^i(\psi,t)$ smoothly depending on $t$ and $\psi$ with the exception of
parallels on ${\Cal F}_{\alpha^*,\beta^*}$ given by $\psi^+_m = \psi^-_{m+1}$.

By $E_i(\psi,t)$ we denote the parallel  transport of $E_i(\psi)$ along $\gamma_\psi$ from an initial point
$\gamma_\psi(0)$ to $\gamma_\psi(t)$. Remind, that above we denoted by $\Phi(\psi)$ the angle of the rotation in
$\Pi(\psi)$ moving $E_i(\psi,0),i=1,2$ to $E_i(\psi,2\pi),i=1,2$. Therefore, the pair of vectors $\tilde
E_i(\psi,t), i=1,2$ obtained from $E_i(\psi,t), i=1,2$ by rotation on angle $t\Phi(\psi)/2\pi$ provides us with
the smooth vector fields over the leaf ${\Cal F}_{\alpha^*,\beta^*}$. In the same way the pair of vectors
$\tilde E_j(\psi,t), j=3,4$ obtained from $E_j(\psi,t), i=3,4$ by rotation on angle $t\Psi(\psi)/2\pi$ provides
us with two more smooth vector fields over the leaf ${\Cal F}_{\alpha^*,\beta^*}$ complementing $\tilde
E_i(\psi,t), i=1,2$ to the (moving) orthonormal base of the tangent space to $M^4$ over the leaf ${\Cal
F}_{\alpha^*,\beta^*}$. Now, $X(\psi,t)$ always denote the unit vector parallel to $\dot\gamma_\psi(t)$,
$Y(\psi,t)$ - unit vector tangent to ${\Cal F}_{\alpha,\beta}$ at the point $\gamma_\psi(t)$ and normal to
$X(\psi,t)$, while $Z(\psi,t)$ and $W(\psi,t)$ - an orthonormal base of the plane normal to ${\Cal
F}_{\alpha,\beta}$ at $\gamma_\psi(t)$. The functions we are going to "smooth" are:
$$
x^i(\alpha^*,\beta^*,\psi,t)=\int\limits_{0}^{t} (\dot\gamma_\psi(t),\tilde E_i(\psi,t)) dt - {{t}\over{2\pi}}\int\limits_{0}^{2\pi} (\dot\gamma_\psi(t),\tilde E_i(\psi,t)) dt, i=1,...,4. \tag 2.1.2
$$
The functions $x^i(\psi,t)$ are $2\pi$-periodic  functions on $\psi$ and $t$. We split them into the Fourier
series:
$$
x^i(\alpha^*,\beta^*,\psi,t)= a^i(\psi)sin(t) + \sum\limits_{k,l=2}^{\infty} a^{i}_k(\psi)sin(kt)+b^{i}_k(\psi)cos(kt), i=1,...,4. \tag 2.1.3
$$
Then our smoothing operators act by cutting out all -  with the exception of the lowest one - harmonics from
these function ("silencing the torus"). Note that if (2.1.3) contains only the first harmonic, or the torus
${\Cal F}_{\alpha^*,\beta^*}$ is "mono-tone", then the corresponding curves $L_\psi(t)$ of the tangent planes in
the Grassmannian fibers are geodesic lines. And vise versa.

We need to check only that this is possible to do by $\delta f$ - by deforming the diffeomorphism $f$
(micro "h"-principle) (and that these deformations diminish our functional ${\Cal H}$). Note that we should consider
only variations of "totally zero charge" - preserving the mean value zero of $x^i$ in (2.1.1). We start with the
following deformations of this kind.

Take the deformation of $f$ localized in a neighborhood of  some point $P$ and given by a $\delta'$-shape vector
field $V=\omega E_i$, where $\omega$ is a function with a small support around $P$, and such that its integrals
over $\gamma_\psi$ are zero. Say, let in coordinates $\{\alpha,\beta,\psi,t\}$ the restriction of $\omega$ on
${\Cal F}_{\alpha^*,\beta^*}$ is: $\omega(\psi,t)=sin(\Lambda (t-t_0)) cos(\Lambda(\psi-\psi_0))$ for some
$0<\psi_0<2\pi, 0<t_0<2\pi$ and $\Lambda$ is sufficiently big such that the support of $\omega$ is inside the
open coordinate rectangular $0<\psi<2\pi, 0<t<2\pi$.

First we note that $E_i(\psi,t)$ are almost stable  under such deformations: if $E_{\delta,i}(\psi,t)$ are the
vector fields defined for $f_\delta$ in the same way as above then the following is true.

\proclaim{Lemma~2.1.1}
$$
\|E_{\delta,i}(\psi,t)-E_{i}(\psi,t)\| \leq \Lambda_1 \delta^2
$$
for some universal constant $\Lambda_1$.
\endproclaim

\demo{Proof} It follows from the Ambrose-Singer formula (1.4): since the "oriented area"  of the deformation
$\gamma_{\delta,\psi}$ is zero - the integral of the restriction of $\omega(\psi,t)$ over $\gamma_\psi$ is zero;
an easy estimate (integration by parts) shows that the norm of the variation of the holonomy operator is small:
$$
\|I_{\gamma_{\delta,\psi}} - I_{\gamma_{\psi}}\| \leq \Lambda_2 supp(\omega) \delta, \tag 2.1.4
$$
where $supp(\omega)$ denotes the diameter of the support of $\omega$, and the  constant $\Lambda_2$ depends on
the norm of the derivative of the curvature operator $\sup
\|\nabla_{\dot\gamma_\delta(t)}R(V_\delta(t),\dot\gamma_\delta(t))\|$. The estimate (2.1.4) implies the claim of
the lemma.

The lemma is proved.
\enddemo

Next we note that from the equality of mixed derivatives
$$
{{\partial^2}\over{\partial\delta \partial t}}f_\delta(\alpha,\beta,\psi,t) = {{\partial^2}\over{\partial t \partial\delta }}f_\delta(\alpha,\beta,\psi,t) \qquad \text{ and } \qquad
{{\partial^2}\over{\partial\delta \partial \psi}}f_\delta(\alpha,\beta,\psi,t) = {{\partial^2}\over{\partial\psi \partial\delta }}f_\delta(\alpha,\beta,\psi,t)
$$
follows that variations of the vectors fields $X,Y$ are given by
$$
\delta X(\psi,t)=\delta \nabla_X V(\psi,t) + O(\delta^2) \qquad \text{ and } \qquad \delta X(\psi,t)=\delta \nabla_X V(\psi,t) + O(\delta^2). \tag 2.1.5
$$
Combining the Lemma~2.1.1 and the last formulas we conclude the following.

\proclaim{Lemma~2.1.2}
$$
|(x_\delta^i(\alpha^*,\beta^*,\psi,t)- x^i(\alpha^*,\beta^*,\psi,t)) - \int\limits_{0}^{t} (\nabla_{\dot\gamma_\psi(t)}V(\psi,t),\tilde E_i(\psi,t)) dt | \leq \Lambda_3 supp(\omega) \delta
$$
for some universal constant $\Lambda_3$.
\endproclaim

\demo{Proof} Indeed, for the proof it is sufficient to note that the variation of the last  term of (2.1.2) (the
mean value of $x^i$) is of the order $\delta supp(\omega)$ by the same "oriented area zero" arguments as in
proving (2.1.4) above.

The lemma is proved.
\enddemo

Now we see that the variation $\delta f$ generated by the vector field
$$
V=\omega E_i, \qquad \text{ where } \qquad \omega(\psi,t)=sin(\Lambda (t-t_0)) cos(\Lambda(\psi-\psi_0)) \tag 2.1.6
$$
gives the following "correction" to the $x^i$:
$$
\delta x^i(\psi,t)= \delta\Lambda sin(\Lambda (t-t_0)) cos(\Lambda(\psi-\psi_0)) + \Lambda_4 {{\delta}\over{\Lambda}}, \tag 2.1.7
$$
for some universal $\Lambda_4$. For further references we may call this "micro-h"-principle.

We defined our variations on the given leaf ${\Cal F}_{\alpha^*,\beta^*}$. In a standard  way it is continued on
its open neighborhood ${\Cal O}_{\alpha^*,\beta^*,\epsilon}$ for an appropriate $\epsilon$ (of the order
$\delta$) preserving estimates above.

Obviously, variations of the type (2.1.6) for large $\Lambda$ generate the space of  all harmonics with
frequencies bigger than $\Lambda$, which would allow us due to (2.1.4) to get rid of all higher harmonics
without increasing much the functional ${\Cal H}$. This would be equivalent to the action of a smoothing
operator $S_\theta$ on $f$.

\subhead 2.2. Basic variations. \endsubhead

For convenience we divide our "corrections" into two types:\footnote{Whether the deformation field $V$ is transversal or tangent to the torus.}

1) "Silencing" as above, i.e., such that $\delta f$ almost does not change the  holonomy operator, see (2.1.4).
In this case the corresponding deformation of the curve of tangent planes $L_\psi(t)$ is the deformation with
(almost) fixed end points;

2) "Tunings" - the new type of deformations changing the holonomy operators (and "shortening" the profile curve
$W(\alpha^*,\beta^*,\psi)$ of the given torus, but almost preserving the integral ${\Cal W(F)}$). In this case
if $L_\psi(t)$ is (almost) a geodesic the corresponding deformation diminishes the distance between its
end-points.

\medskip

\subhead 2.2.1 "Silencing". Deformations straightening the curve of tangent planes \endsubhead

Two bundles over $M^4$: of the Grassmannians $G$ and isometries $SO$ are associated to the tangent bundle $TM^4$
of $M^4$ and naturally inherit the connection, or the parallel translation which we denote again by
$I_\gamma(t)$. The parallel translation identifies restrictions $G(t)$ and $SO(t)$ of the bundles $G$ and $SO$
to the curve $\gamma$ with the point $P$ deleted with the direct metric products $G_P\times (0,2\pi)$ and
$SO_P\times (0,2\pi)$; or with the metric products with closed interval, if we double the point $P$ and fibers
over it. Correspondingly the section $T(t)$ of the Grassmann bundle $G(t)$ of two-planes over $\gamma$ is
identified with the section $(L(t),t)$ of the direct product $G_P\times [0,2\pi]$. Slightly abusing notations,
below we use this identification with the direct metric products when working in $G(t)$ and $SO(t)$. E.g.,
covariant derivative $\nabla_{\dot\gamma(t)}$ now is the derivative on $t$ in the metric products $G_P\times
[0,2\pi]$ and $SO_P\times [0,2\pi]$, while $\{X(t),Y(t)\}$ is an orthonormal base of $L(t)$; $L(t)=X(t)\wedge
Y(t)$.

\medskip

The  tangent space to the Grassmann fiber $G(t)$ at the point $X(t)\wedge Y(t)$ is generated by four vectors
$A_{XZ}(t)=X(t)\wedge Y(t)$, $A_{XW}(t)=X(t)\wedge W(t)$, $A_{YZ}(t)=Y(t)\wedge Z(t)$ and $A_{YW}(t)=Y(t)\wedge
W(t)$. These vectors are images under the natural projection $\pi: SO(t)\to G(t)$ of the generators of isometry
transformations of $T_\gamma(t)M^4$ - rotations in corresponding planes, e.g., let $SO_{XZ}(\phi)$ be a rotation
with a unit speed of $T_\gamma(t)M^4$ in the plane $X(t)\wedge Z(t)$ (about the fixed plane $Y(t)\wedge W(t)$).
Then the vector $A_{XZ}(t)$ equals the image of the derivative of $SO(\phi)$ at $\phi=0$ under the differential
of $\pi$. Naturally, we assume that $G$ is endowed with the metric such that $\pi: SO(t)\to G(t)$ is a
Riemannian submersion where the metric on $SO$ is given by the standard Lipschitz-Killing form: for a skew
symmetric matrix $A=(A_{ij})$ it holds $\|A\|^2=\sum_{i<j}A^2_{ij}$. In this metric
$\{A_{XZ},A_{XW},A_{YZ},A_{YW}\}$ is an orthonormal base of the tangent space to the Grassmannian fiber $G$.

For a  given instant $t$ the derivative of $L(t)$ on $t$ may be written as $S(t)L(t)$, where $S(t)$ is the
generator of a one-parameter subgroup $SO(t,\epsilon)$ of $SO(t)$ such that
$$
L(t+\epsilon)=SO(t,\epsilon)L(t) + o(\epsilon) \qquad \text{ where } \qquad S(t)={{d}\over{d\epsilon}}SO(t,\epsilon)_{|\epsilon =0} \tag 2.2.1
$$
($S(t)$ stands for the "second form" of $\Cal F$ in the direction $\dot\gamma(t)$.)

We  denote by $\pi(t)$ and $n(t)$ invariant subspaces of $SO(t,\epsilon)$ and by $\{e_i\}, i=1,2,3,4$ an orthonormal base of $T_{\gamma(t)}M$ such that $e_i, i=1,2$ generates $\pi(t)$, while $e_i,i=3,4$ - $n(t)$. In this base $SO(t,\epsilon)$ is given by block-diagonal matrix acting (with respect to $\epsilon$) on invariant subspaces by rotations on $\epsilon\phi(t)$ and $\epsilon\psi(t)$ correspondingly; where $\phi(t)$ and $\phi(t)$
are some constants (speed of rotations) depending on $t$:
$$
SO(t,\epsilon) =
\left(
\matrix
\cos(\epsilon\phi(t)) & -\sin(\epsilon\phi(t)) & 0& 0\\
\sin(\epsilon\phi(t)) & \cos(\epsilon\phi(t)) & 0& 0\\
0& 0& \cos(\epsilon\psi(t)) & -\sin(\epsilon\psi(t)) \\
0& 0& \cos(\epsilon\phi(t)) & -\sin(\epsilon\phi(t))\\
\endmatrix
\right)
\tag 2.2.2
$$
To simplify computations instead of $\{A_{XZ},A_{XW},A_{YZ},A_{YW}\}$ we use  another orthonormal base
$e_{ij}=e_i\wedge e_j, i,j=1,2,3,4$ provided by $e_i$. In this base the geodesic curvature of the curve $L(t)$
is easy to compute.

\proclaim{Lemma~2.2.1} For the geodesic curvature $Gc(t)$ of the curve $L(t)$ it holds
$$
Gc(t)=(\phi^2(t)+\psi^2(t)) \|R(t)\|, \tag 2.2.3
$$
where $R(t)=R_{ij}(t), 1\leq i\leq 2<j\leq 4$ is the matrix of (covariant) derivatives
$$
R_{ij}(t)=(({{d}\over{dt}}S(t))e_i,e_j). \tag 2.2.4
$$
\endproclaim

\demo{Proof}  Geodesic curvature of the curve $L(\epsilon)$ measures its deviation from the geodesic line: if
$\bar L(\epsilon)$ is the geodesic issuing from the same initial point $L(0)$ and speed $\dot L(o)$ as the curve
$L$, then the vector $L(\epsilon)\bar L(\epsilon)$ (i.e., the unit direction of the minimal geodesic going from
$L(\epsilon)$ to $\bar L(\epsilon)$ multiplied by the distance between $L(\epsilon)$ and $\bar L(\epsilon)$)
equals (up to higher order terms) the vector $Gc$ of the geodesic curvature of $L$ at $L(0)$ multiplied by
$\epsilon^2/2$ times the square of the speed $\|\dot L(0)\|^2$. To see this - take the second derivative of
$L(\epsilon)\bar L(\epsilon)$ on $\epsilon$ at zero. For instance, what is the geodesic curvature of the curve
of isometries $SO(\tau)$ sending $L(0)$ to $L(\tau)$ at the instance $\tau=t$? From (2.2.1) it follows that the
geodesic we should compare to is the one-parameter family $SO(t,\epsilon)$ of isometries, the speed equals
$\|S(t)\|=(\phi^2(t)+\psi^2(t))$; and the second derivative of the vector $SO(t+\epsilon)SO(t,\epsilon)$ on
$\epsilon$ at zero is ${{d}\over{dt}}S(t)$. The last vector - tangent to the bundle $SO$ of isometries - under
the submersion $\pi$ to $G$ goes to the vector $R(t)$, from which the claim of the lemma follows.

Lemma~2.2.1 is proved.
\enddemo

\medskip

The idea is to  construct the variation field $V$ of $f_\delta$ such that the corresponding variation $\delta
L(t)$ is the vector parallel to the vector $Gc(t)$ of the geodesic curvature of the curve $L(t)$. Then the first
variation formula would imply that the length of $L(t)$ is decreasing under variation $f_\delta$. An arbitrary
vector tangent to $G(t_0)$ at the point $T(t_0)=X(t_0)\wedge Y(t_0)$ is a linear combination of four vectors
$A_{ZY}(t_0)$, $A_{WY}(t_0)$, $A_{XZ}(t_0)$ and $A_{XW}(t_0)$. Denote by $V_{ZY}(t_0)$ the variation field
$$
V_{ZY}(t_0)(\psi,t)={{1}\over{\Lambda}}\sum\limits_i sin(\Lambda(t-t_0))cos(\Lambda(\psi-\psi_0)) (Z(\psi_0,t_0),E^i(\psi_0,t_0)) E^i(\psi,t) \tag 2.2.5
$$
for coordinates  $(\psi,t)$ in the coordinate rectangular $\Lambda |t-t_0|\leq\pi$ and $\Lambda
|\psi-\psi_0|\leq\pi/2$, and zero elsewhere. From (2.1.5) it follows that over the same rectangular
$$
{{{\partial}\over{\partial\delta}}X(\psi,t)_\delta }_{| \delta=0} = cos(\Lambda(t-t_0))cos(\Lambda(\psi-\psi_0)) Z(t) + O(\delta)
$$
and
$$
{{{\partial}\over{\partial\delta}}T(\psi,t)_\delta }_{| \delta=0} = cos(\Lambda(t-t_0))cos(\Lambda(\psi-\psi_0)) Z(t)\wedge Y(t) \tag 2.2.6
$$
correspondingly for sufficiently small $\delta$ or large $\Lambda$ ($\delta\Lambda <1$). In order to have
smooth $f_\delta$ diffeomorphisms we smooth a little bit more $V$ in, say, $1/10 \Lambda$-neighborhood of the
coordinate rectangular and obtain the variation we need. In the same way we construct the variational field
$V_{WY}$ giving the variation of $T(t)$ in the direction $W\wedge Y$; and in a little bit more subtle way -
because $Y$ is not the $\psi$-coordinate field, and interchanging functions $sin(t)$ and $cos(\psi)$ with
$cos(t)$ and $sin(\psi)$ if necessary; variational fields $V_{XZ}$ and $V_{XW}$ such that for the corresponding
variations of the tangent plane $T(t)$ are parallel to $A_{XZ}$ and $A_{XW}$ correspondingly. We summarize our
arguments as follows.

\proclaim{Lemma~2.2.2} For  an arbitrary vector $R_{ij}(t_0)$ tangent to $G(t)$ at $T(t)$ there exists a tame
deformation of $f$ given by some tame vector field $V_{ij}$ such that
$$
{{{\partial}\over{\partial\delta}}T(\psi,t)_\delta }_{| \delta=0} = \omega(\psi,t) R_{ij}(t)
$$
where $\omega(\psi,t)$ - "good" bump function with vanishing integrals over all $\gamma_\psi$.
\endproclaim

Using partition of unity an  arbitrary vector field $R(t)$ tangent to $G(t)$ along $T(t)$ may be presented as
the sum of vector fields with small supports. Denote by $\omega$ the maximum size of such supports, then
construct $R(t)$ (sufficiently smooth) approximating the unit vector field parallel to the geodesic curvature
vector $Geod(t)$ of $T(t)$. Then the arguments above by linearity of the derivative lead to the following
estimates:
$$
{{\partial}\over{\partial\delta}} {\Cal H}(\gamma_{\delta,\psi}(t)) <  -(1-\Lambda_5)\int\limits_0^{2\pi} \|Geod(t)\| dt + \Lambda_5 \omega. \tag 2.2.7
$$
for some universal positive  $\Lambda_5$. Since we may choose the parameter $\omega$ converging to zero during
the corrections the last estimate will guarantee the convergence of the Nash process.

The constructed  corrections diminish ${\Cal H}(\gamma_{\psi_0}(t))$ of the one particular curve
$\gamma_{\psi_0}$. By continuity the same is true for $\gamma_{\psi}$ for $\psi$ from some $\omega'$-neighborhood of $\psi_0$; i.e., if the deformation vector field $V$ implies the deformation $\delta T(\psi_0,t)$ with an integral positive scalar product with a geodesic curvature vector of $T(\psi_0,t)$, then the scalar product of $\delta T(\psi,t)$ with the geodesic curvature vector of $T(\psi,t)$ is still positive for $\psi$ from some  $\omega'$-neighborhood of $\psi_0$. Taking $\delta$ so small that the support of the deformation $\delta f$ lies in this neighborhood we define the deformation which diminish ${\Cal H}$ in general case, when the extremal $\gamma^*$ is not necessarily isolated. Remark also that in the process of diminishing ${\Cal H}(\gamma^*)$ we may at the same time deform $f$ along other $\gamma_\psi(t)$ diminishing ${\Cal H}(\gamma_\psi)$. The limit of such deformations is the torus ${\Cal F}^*$ where all curves $L_\psi(t)$ of tangent planes are geodesics. We call such torus ${\Cal F}^*$ --- "mono-tone".

Our next step is to deform this "mono-tone" torus into the "silent" one.

\medskip

\subhead 2.2.2 "Tunings". Deformations of the mono-tone torus \endsubhead

The sequence of corrections  silencing ${\Cal F}_{\alpha,\beta}$ gives us the torus such that all (corresponding
to extremals) curves $T_\psi(t)$ of tangent planes are geodesics. In this case, whether the variation
$T_\delta(t)$ of a geodesic line shortens it or not - depends on the difference of scalar products of the
vectors of variations $\partial T_\delta/\partial\delta$ with the unit speed $\dot T(t)$ of $T$ at its end
points. Therefore, to further shorten $T(t)$ and diminish ${\Cal H}$ we may do the following.

{\bf "Boundary operator on flag foliations".} Take a meridian $l(\psi)=\gamma_\psi(0)$ of a torus and consider the
Grassmann bundle of two-planes over it. This bundle has two sections $\Pi(\psi)$, $N(\psi)$ we already defined;
and two sections $T_\psi(0)$ and $T_\psi(2\pi)=I_{\gamma_\psi}(T_\psi(0))$ - of end points of the curves of
tangent planes. For each $\psi$ in the fiber $G_\psi$ of this bundle the parallel translation $I_{\gamma_\psi}$
acts as a rotation about $\Pi(\psi)$ and $N(\psi)$ on angles $\Phi(\psi)$ and $\Psi(\psi)$ correspondingly,
i.e., for some one-parameter group of rotations $SO_\psi(t)$ it holds $SO_\psi(2\pi)=I_{\gamma_\psi}$; while
$T_\psi(t)$ is an $0\leq t\leq 2\pi$-interval of the orbit of the action of this group on $G_\psi$. We may
further simplify the picture by doubling the point $l(0)=l(2\pi)$, identifying the Grassmann bundle $G_\psi$
over the closed interval $0\leq \psi\leq 2\pi$ with the direct product $[0,2\pi]\times G$ in such a way that
$\Pi(\psi)$ and $N(\psi)$ be two parallel sections and projecting to $G$. Then in the Grassmannian $G$ we will
have the curve $L^-(\psi)$ of tangent planes to the torus along the meridian $l(\psi)$ - the beginning points of
$T_\psi(t)$, while the curve of end points $L^+(\psi)=T_\psi(2\pi)$ is obtained from $L^-(\psi)$ by the rotation
$SO_\psi(2\pi)$; i.e., these end points are connected by the geodesic $T_\psi(t)=SO_\psi(t) L^{-}(\psi), 0\leq
t\leq 2\pi$.

Note that the  curve $L^-(\psi)$ in $G$ is closed: $L^-(0)=L^-(2\pi)$. Indeed, $\Pi(\psi)$ and $N(\psi)$ are
generated by $\tilde E_i(\psi,t)$ for $t=0$ and coincide for $\psi=0,2\pi$. This identifies $\{0\}\times G$
with $\{2\pi\}\times G$ in the direct product bundle $[0,2\pi]\times G$ with $\Pi(\psi)$ and $N(\psi)$ parallel
sections.

We may clarify this picture further. Consider the direct product $[0,2\pi]\times T_{P}M^4$ as the trivial bundle over
meridian $l(\psi)$ with the point $P=l(0)$ doubled, and vector fields $\tilde E_i(\psi,0), i=1,...,4$ parallel.
Note also that our condition ${\Cal W}(\alpha,\beta)=0$ implies that the vector fields $\tilde E_i(\psi,2\pi)$
in such a bundle have equal initial and end values: $\tilde E_i(0,2\pi)=\tilde E_i(2\pi,2\pi)$ since the angle
between $\tilde E_i(\psi,0)$ and $\tilde E_i(\psi,2\pi)$ in the corresponding invariant plane equals
$\Phi(\psi)$ or $\Psi(\psi)$ who's integrals on $\psi$ vanish.

On the geometry of $G$: Remind that $G$ can be identified with the direct product $S^2\times S^2$ as
follows: Let $\Pi$ is generated by $E_1,E_2$, $N$ by $E_3,E_4$. Then any particular point $L$ in $G$ has the
basis $e_1,e_2$ such that $e_1=E_1 cos(\lambda)+E_3 sin(\lambda)$ and $e_2= - E_2 sin(\mu)+E_4 cos(\mu)$, and
$L$ is the image of $\Pi$ under the composition of rotations on $\lambda$ about the plane $\{E_2,E_4\}$ and on
$\mu$ about the plane $\{E_1,E_3\}$. For a given $(\lambda,\mu)$ the set of corresponding $L$ is the torus
$T^2_{\lambda,\mu}$ in $G$ (similar to the above $\bar{\Cal F}_{\lambda,\mu}$), and $\sqrt{\lambda^2+\mu^2}$
gives the distance from $\Pi$ to $L$. The distance from $L$ to $N$ is $\sqrt{(\pi-\lambda)^2+(\pi-\mu)^2}$, so
that $T^2_{\pi/2,\pi/2}$ is the equidistant set of $\Pi$ and $N$.

Consider these $\lambda(\psi),\mu(\psi)$ functions: they have equal end-values:
$$
\lambda(0)=\lambda(2\pi)\qquad \text{ and } \qquad \mu(0)=\mu(2\pi)
$$
and, considering the functional
$$
\int\limits_0^{2\pi} (\lambda'(\psi))^2+(\mu'(\psi))^2 d\psi \tag 2.2.8
$$
we may diminish it by construct deformations $\delta f$  as above. This results in the mono-tone torus $\Cal F^*$ which is "tuned";
i.e., such that the "left-end" curve $L^-(\psi)$ would belong to some particular torus $T^2_{\lambda,\mu}$. The parallel translations then
draw a geodesic $L(t)$ inside this torus connecting $L^-(\psi)$ with $L^+(\psi)$.

Another thing is that we may assume that this torus is $T^2_{\pi/2,\pi/2}$, i.e., is equidistant to $\Pi$ and $N$.

Indeed. It is easy to see that, for the  mono-tone torus its ${\Cal H}(\psi)$ functional is determined by the "picture"
consisting of data $L^-(\psi)$ and $SO_\psi$:
$$
{\Cal H}(\psi)=dist_G(L^-(\psi), SO_\psi(L^-(\psi))). \tag 2.2.9
$$

Next we note that in order to diminish $dist_G(L^-(\psi), SO_\psi(L^-(\psi))=L^+(\psi))$ we may construct the
deformations $\delta f$ which moves the curve $L^-(\psi)$ to one of the fixed points $\Pi$ or $N$ of all
$SO_\psi$ rotations on $G$. To be more precise, first we check that if $dist_G(\Pi,L^-(\psi))\not=
dist_G(N,L^-(\psi))$, then the deformation $\delta f$ such that again: 1) the holonomy is almost stable (with zero
"oriented area" as above), and 2) $L^-_\delta(\psi)$ moves toward $\Pi$ if $dist_G(\Pi,L(\psi)) <
dist_G(N,L(\psi))$, or toward $N$ if $dist_G(\Pi,L(\psi)) > dist_G(N,L(\psi))$ - diminish ${\Cal H}(\psi)$,
i.e., the distance between end points $T_\psi(0)$ and $T_\psi(2\pi)$.

Indeed, connect these  endpoints $L^-(\psi)$ and $L^+(\psi)$ with $\Pi$ by some minimal geodesics $L^-(\psi,s)$
and $L^+(\psi,s)$, and consider a deformation $f_\delta$ such that the corresponding $\delta L^-(\psi)$ equals
the unit direction of $\dot L^-(\psi,s)$ at $s=0$. If, again, the variation of the holonomy operator
$I_{\gamma_\psi}$ under such deformation vanishes (is "small"), then $\delta L^+(\psi)$ equals the unit
direction of $\dot L^+(\psi,s)$ at $s=0$; i.e., under such deformation the end points $L^-(\psi)$ and
$L^+(\psi)$ moves along the sides of the triangle in $G$ with vertices $\Pi$, $L^-(\psi)$ and $L^+(\psi)$. When
$dist_G(\Pi,L(\psi)) < dist_G(N,L(\psi))$ the angles $L^-(\psi)$ and $L^+(\psi)$ of this triangle are less than
$\pi/2$, and the first variation formula implies
$$
dist_G(L^-_\delta(\psi), L^+_\delta(\psi)) < dist_G(L^-(\psi), L^+(\psi)). \tag 2.2.10
$$
Is it possible to construct such a deformation, i.e., satisfying two properties above? Yes, such deformation
$\delta f$ may be constructed with a small support (take $\|V(t)\|=sin(t/\epsilon)$ on
$-\pi\epsilon<t<\pi\epsilon$ and zero elsewhere) the change of the holonomy operator is small $(\| \delta I
\|\leq O(\delta\epsilon))$. If the angles $T(0)$, $T(2\pi)$ of the triangle $\Delta (\Pi T(0) T(2\pi))$ in $G_P$
(which are equal) are strictly less than $\pi/2$ this gives (2.2.10) for $\epsilon$ sufficiently small.\footnote{When those
angles  equal $\pi/2$ the first derivative of the distance in $G_P$ under the deformation above is
zero, but the second is negative. Therefore, (2.2.10) would follow if we could further (up to
$O(\delta^2\epsilon)$) reduce the affect of the change of the holonomy. Indeed, assuming  that under some local
deformation $\delta f$ the "extra-push" in the holonomy operator generated by $R(\delta f,X)$, as we described
above, moves $T(2\pi)$ not normal to the direction of $\dot T(t)$ at $T(2\pi)$, then, changing $\delta f$ to
$-\delta f$ if needed, we obtain shortening deformation. If not, i.e., all $\delta T(2\pi)$ under $\delta
I_\gamma$ are normal to $\dot T(2\pi)$; then the change of the holonomy has reduced affect. Still, to obstacle such deformation
might the impossibility to choose correctly the deformation direction toward $\Pi$ or $N$ for all $\psi$ together.}

We formulate the obtained result as follows.

\proclaim{Lemma~2.2.3} If $dist_G(\Pi,L(\psi^*)) \not= dist_G(N,L(\psi^*))$ then there always exist a
correction $\delta f$ diminishing the ${\Cal H}(\psi)$ for all $\psi$ from some neighborhood of $\psi^*$;
i.e., such that (2.2.10) holds. Therefore, (from now on till the end of this section) we may assume that
the angles $L^-(\psi)$ and $L^+(\psi)$ in the triangle $\triangle{\Pi L^-(\psi)L^+(\psi)}$ in $G$ equal
$\pi/2$, or that $L(\psi,t)$ lie in the torus $T^^2_{\pi/2,\pi/2}$ in $G$ equidistant to $\Pi$ and $N$.
\endproclaim

2.4. The final steps include the "fine-tuning" deformations along the tangent vector fields $X$ and $Y$; i.e., which is equivalent to
re-parametrization of the extremal torus ${\Cal F}^*$. To picture these tunings, call the union $A(\psi,t)$ of all curves $L_\psi(t)$ (which are geodesics $L_\psi(t)=SO_\psi(t)L^-(\psi)$ in $G$) connecting $L^-(\psi)$ with $L^+(\psi)$ - the annulus of ${\Cal F}^*$ (on $T^2_{\pi/2,\pi/2}$ in $G$) with boundary curves $L^-(\psi)$ and $L^+(\psi)$. Then deformations along $X$ will change the inner boundary $L^-(\psi)$ to $L_\psi(t+tune(\psi))$ with the corresponding move of the outer boundary: $L^+(\psi)$ to $L_\psi(2\pi+tune(\psi))$, where $L_\psi(t)$ for $t>2\pi$ is the point on the geodesic $SO_\psi(t)L^-(\psi)$ which is defined for all $t$. The deformations along $Y$ change the parameter $\psi\to\psi_\delta$ - which will result in the different correspondence between start and end points: under such deformation the end  $L_\delta^+(\psi)$ is $L^+(\psi_\delta)$, see details below.

By definition ${\Cal H}(\psi)$ equals the distance from the start point $L^-(\psi)$ to the end-point $L^+(\psi)$ in $A(\psi,t)$. Now, a simple picture of the deformations of $A(\psi,t)$ as above: 1) moves along geodesics $L_\psi(t)$, and 2) independent "rotations" of boundary curves $L^-(\psi)$ to the end-point $L^+(\psi)$ in $A(\psi,t)$ by changing the parameter $\psi$ on them; proves that ${\Cal H}(\psi)$ could be easily diminished unless 1) it is identically zero, or 2) the inner and outer curves $L^-(\psi)$ and $L^+(\psi)$ are two points connected by a single geodesic $L(t)=L_\psi(t)$ which implies that the vector ${\Cal W}(\psi)$ is constant and not zero. Since the last case contradicts our condition (2.2.11) (integral of ${\Cal W}(\psi)$ on $\psi$ is zero) this proves that our "fine-tuning" will ultimately results in the "silent" torus ${\Cal F}^*$ with ${\Cal H(F)}\equiv 0$.
This proves the Lemma~2.2.4 below.

Let us provide some details for this procedure:

How the deformation may  affect the holonomy operator? Again, by Ambrose-Singer formula the local deformation
$\delta f$ around the point $\gamma(t_0)$ leads to the variation of $\gamma$ and changes the holonomy by an
"extra-push" $\delta I_\gamma=exp(R(\delta V, X))$ "inserted" in $\gamma(t_0)$. For instance, by definition
$T(2\pi)$ is the parallel transport of $T(0)$ along $\gamma$, and after the deformation $\delta f$ will be
generated by the vectors $X+\delta X$ and $Y+\delta Y$, where $\delta X$ equals the integral of $R(\delta V,
X)X$, and $\delta Y$ equals the integral of $R(\delta V, X)Y$. If the restriction of the deformation vector
field $V$ on $\gamma(t)$ is supported on the interval $(t_0-\epsilon,t_0+\epsilon)$ then all $T(t)$ and $L(t)$
for every $t$ from the interval $[t_0-\epsilon, \pi]$ may be affected; i.e., the local variation $\delta f$ may
lead to the global deformation of the curve $T(t)$.

Assuming that the curve  $T(t)$ after all our deformations above is already a geodesic the first variation
formula says that the deformation $\delta f$ shortens its length (i.e., (2.10) holds) if the direction
$R(\delta f, X)$ of  the generator of the additional rotation $\delta I$ has negative scalar
product with the direction $\dot T(2\pi)$. \footnote{Correspondingly, the stationary $\gamma^*$ of the
functional ${\Cal H}$ provides the orthogonality equalities for the curvature tensor of $M^4$.}

\medskip

Remind, that for a given torus
${\Cal F}_{\alpha,\beta}$ its profile curve is ${\Cal W}((\alpha,\beta,\psi))=(\Phi(\alpha,\beta,\psi),\Psi(\alpha,\beta,\psi))$
where each parallel translation $I_{\gamma_\psi}$ acts on its invariant planes $\Pi(\psi)$ and $N(\psi)$ by rotations on angles
$\Phi(\alpha,\beta,\psi)$ and $\Psi(\alpha,\beta,\psi)$ correspondingly. Remind also that we consider only the torus such that
$$
{\Cal W}(\alpha,\beta)=\int\limits_0^{2\pi}(\Phi(\alpha,\beta,\psi),\Psi(\alpha,\beta,\psi)) d\psi=0, \tag 2.2.11
$$
or with a profile which center of gravity lies in the origin. Next we omit $\alpha,\beta$ from notations for simplicity and denote
by $w(\psi)$ the speed $w(\psi)={\Cal W}'(\psi)$ of the profile curve in the plane $\{(\Phi,\Psi)\}$.

Reprove the Ambrose-Singer formula:

Let $\gamma_\delta(t)$  be some deformation of the curve $\gamma(t)$ with  $V(t)={\partial
(\gamma_\delta(t))/\partial\delta}_{|\delta=0}$ - the deformation field. Let $E(\delta,t)$ be a vector field
parallel along $\gamma_\delta(t)$; i.e.,
$$
\nabla_{\partial/\partial t}\tilde E(\delta,t)\equiv 0.
$$
Since vector fields  $\partial/\partial t$ and $\partial/\partial \delta$ commute it holds by the definition of
the curvature tensor:
$$
\nabla_{\partial/\partial t}\nabla_{\partial/\partial \delta} E(\delta,t)=R({\partial/\partial t},{\partial/\partial \delta})E(\delta,t).
$$
If $D(t)$ is any vector  field along $\gamma(t)$ and $\bar D(t)$ denotes its parallel transport $I_\gamma(t)$
from the point $\gamma(t)$ to $\gamma(0)$ along $\gamma$, then, again by the definition of the parallel
transport it holds
$$
{{d}\over{dt}} \bar D(t)= I_\gamma(t) \nabla_{\partial/\partial t} D(t)\qquad \text{ and }\qquad \bar D(2\pi)-D(0)= \int\limits_0^{2\pi} I_\gamma(t)(\nabla_{\partial/\partial t} D(t) dt.
$$
Applying this to the  vector field $D(t)=\nabla_{\partial/\partial \delta} E(\delta,t)_{\delta=0}$ we get the
Ambrose-Singer formula (1.4) as follows:
$$
I_\gamma(2\pi) (\nabla_{\partial/\partial \delta} E(\delta,2\pi)_{\delta=0}) -  \nabla_{\partial/\partial \delta} E(\delta,0)_{\delta=0}= \int\limits_0^{2\pi} I_\gamma(t) (R({\partial/\partial t},{\partial/\partial \delta})E(0,t)) dt. \tag {2.2.12 Ambrose-Singer}
$$

\medskip

Now we apply this  formula to the vector fields $\tilde E_i(\psi,t)$ and obtain formulas for the speed of the profile curve.
Recall that these fields are obtained by rotating with a constant speed   $\Phi(\psi)/2\pi$ and $\Psi(\psi)/2\pi$ the pairs $E_i(\psi,t)$
for $i=1,2$ and $3,4$ (bases of invariant planes $\Pi$ and $N$) correspondingly of parallel vector fields along $\gamma_\psi(t),
0\leq t\leq 2\pi$ in order to have smooth vector fields on closed $\gamma_\psi(t),0\leq t\leq 2\pi$ with
$\gamma_\psi(0)=\gamma_\psi(2\pi)$. By this definition $\tilde E_i(\psi,t), i=1,2$ and $\tilde E_i(\psi,t),
i=3,4$ are projections on invariant planes $\Pi$ and $N$ correspondingly of orbits under the action of the one
parameter group $\tilde SO(\psi,t)=exp(t\tilde S(\psi)/2\pi)$ generated by the block matrix $\tilde S(\psi)$,
which in the basis $E_i$ has the form:
$$
\tilde S(\psi)=\Phi(\psi) \left( \matrix 0& -1 & 0 & 0\\
1&0&0&0\\
0&0&0&0\\
0&0&0&0\\
\endmatrix \right) +
\Psi(\psi)\left( \matrix 0& 0 & 0 & 0\\
0&0&0&0\\
0&0&0&-1\\
0&0&1&0\\
\endmatrix \right)
$$
The Ambrose-Singer  formula above gives the variation of $\tilde SO(\psi,t)$, and applying the first variation
formula we conclude that the derivatives of $\Phi,\Psi$ are the scalar products of the derivative of the
holonomy operator $I_{\gamma_\psi}$ on $\psi$ with this vector: i.e., the following is true:
$$
\Phi'(\psi)= \int\limits_0^{2\pi} I_\gamma(t) (R({\partial/\partial t},{\partial/\partial \psi})E_1(\psi,t),E_2(\psi,t)) dt
$$
and
$$
\Psi'(\psi)= \int\limits_0^{2\pi} I_\gamma(t) (R({\partial/\partial t},{\partial/\partial \psi})E_3(\psi,t),E_4(\psi,t)) dt
$$
If we denote by  $\sigma(\psi,t)$ the area element of the torus ${\Cal F}$ in coordinates $\{\psi,t\}$, and by
${\Cal R}_\Pi(\psi,t)=I_\gamma(t) (R({\partial/\partial t},{\partial/\partial \psi})E_1(\psi,t),E_2(\psi,t))$
and ${\Cal R}_N(\psi,t)=I_\gamma(t) (R({\partial/\partial t},{\partial/\partial \psi})E_3(\psi,t),E_4(\psi,t))$
we get the following formula:
$$
\Phi'(\psi)= \int\limits_0^{2\pi} \sigma(\psi,t) {\Cal R}_\Pi(\psi,t) dt\tag 2.2.13
$$
and
$$
\Psi'(\psi)= \int\limits_0^{2\pi} \sigma(\psi,t) {\Cal R}_N(\psi,t) dt,\tag 2.2.14
$$
where ${\Cal R}_\Pi(\psi,t)$  and ${\Cal R}_N(\psi,t)$ are functions on the torus and do not depend on its particular parametrization by $\{\psi,t\}$
like $\sigma(\psi,t)$. Of course, the same formulas  are true for an arbitrary deformation $\delta f$ given by the vector field $V$, i.e.,
$$
{{\partial\Phi_\delta}\over{\partial\delta}}(\psi)_{|\delta=0}= \int\limits_0^{2\pi} I_\gamma(t) (R({\partial/\partial t},V(\psi,t))E_1(\psi,t),E_2(\psi,t)) dt
\tag 2.2.15
$$
and
$$
{{\partial\Psi_\delta}\over{\partial\delta}}(\psi)_{|\delta=0}= \int\limits_0^{2\pi} I_\gamma(t) (R({\partial/\partial t},V(\psi,t))E_3(\psi,t),E_4(\psi,t)) dt.
\tag 2.2.16
$$

\medskip

"Fine-Tunings" along $X$. What happens under deformations along $X$ is obvious: if the vector field $V$ is parallel to $X$, i.e., $V(\psi,t)=tune(\psi,t)X$; then $\gamma_\delta(\psi,0)=\gamma(\psi, \delta tune(\psi,0))$. At this point the invariant planes $\Pi(\psi,\delta tune(\psi,0))$ and $N(\psi,\delta tune(\psi,0))$ are parallel transports of $\Pi(\psi,0)$ and $N(\psi,0)$ along $\gamma_\psi$, while the tangent plane $T_\delta(\psi,0)=T(\psi,\delta tune(\psi,0))$ is the image of the parallel transport of $T(\psi,0)$ under rotation $SO_\psi(\delta tune(\psi,0))$, i.e., corresponds to $L(\psi,\delta tune(\psi,0))$ in our trivialization of the Grassmann bundle as above. The corresponding move of the outer boundary is: $L^+(\psi)$ to $L_\psi(2\pi+tune(\psi))$, where $L_\psi(t)$ for $t>2\pi$ is the point on the geodesic $SO_\psi(t)L^-(\psi)$ which is defined for all $t$.

In all we see that considered deformations $\delta f$ along $X$ move inner and outer curves of $A(\psi,t)$ in the directions of geodesics $L(\psi,t)$ at $t=0,2\pi$ correspondingly in $G$ with a speed given by the functions $tune(\psi,0)$ and $tune(\psi,2\pi)$.

"Fine-Tunings" along $Y$ we define as $\delta f$ generated by
$$
V(\psi,t)= tune(\psi,t){{\partial}\over{\partial \psi}}(\psi,t),  \tag 2.17
$$
with $tune(\psi,t)$ smooth equals one (or minus one) over all $\psi^*-\epsilon<\psi<\psi^*+\epsilon$ and $\epsilon<t<2\pi-\epsilon$ but vanishing on the boundary of the coordinate rectangular $[\psi^*-2\epsilon, \psi^*+2\epsilon]\times [0,2\pi]$ for small $\epsilon$.
It holds
$$
{{\partial\Phi_\delta}\over{\partial\delta}}(\psi)_{|\delta=0}= \int\limits_0^{2\pi}
tune(\psi,t) \sigma(\psi,t) {\Cal R}_\Pi(\psi,t) dt
\tag 2.2.18
$$
and
$$
{{\partial\Psi_\delta}\over{\partial\delta}}(\psi)_{|\delta=0}= \int\limits_0^{2\pi}
tune(\psi,t)\sigma(\psi,t) {\Cal R}_N(\psi,t) dt.
\tag 2.2.19
$$
Comparing with (2.2.13) and (2.2.14) we see that the effect of such tuning is that we act on the curve $L^-(\psi)$ not by rotations $SO(\psi,t)$ but almost by $SO(\psi+\epsilon,t)$ (or $SO(\psi-\epsilon,t)$), resulting effectively in the change of the parameter $\psi$ to $\psi+\epsilon$ (or $\psi-\epsilon$ correspondingly) on the inner boundary curve. If the angle between $L^-(\psi)$ and the direction of $L_\psi(t)$ at $t=0$ is acute such tuning is diminishing the distance between end-points of $L_\psi(t)$ (if the angle is obtuse - use the tuning almost equal to minus one).

Taking three points $L^-(\psi_i), i=1,2,3$ not on the same line we easily balance tunings defined in neighborhoods of these points to keep the center of gravity of the profile curve at zero, and shortening profile at the same time.

Again, Obviously, the only  stationary $\omega$ under described deformations are straight segments. But if $\omega$ is closed curve, as in our case, the deformation is always possible.

The "profile straightening"  is a tame deformation. Indeed: $\delta f$ might be very big compare to the change of the profile if only the curvature generators vanish. But then the curve $\omega$ would be of almost zero curvature. On the other hand - it is inside some disk in the plane - in its vertices the curvature is less than the inverse radius of the boundary circle - which is bounded.

We formulate the obtained result as follows.

\proclaim{Lemma~2.2.4. (On "Local Minimum")} If ${\Cal H(F^*)}\not= 0$ for the mono-tone torus ${\Cal F}^*_{\alpha,\beta}$ such that ${\Cal W}(\alpha,\beta)=0$, then there exists the ("fine-tuning") deformation $\delta f$ diminishing ${\Cal H(F^*)}$.
\endproclaim

\medskip

Let us repeat again: it is not clear immediately that our deformations are smoothing operators as we defined in the section 2.1 above.
For our convenience we defined them in $\{\alpha,\beta,\psi,t\}$-coordinates, i.e., which is different from the definition in 2.1.
As we note above, we should consider the isometric embedding $M^4\to R^N$ and "translate" our constructions to Euclidean coordinates in $R^N$.

\medskip

\head 3. The Nash process \endhead

\subhead 3.0. Standard assumptions \endsubhead

According to our general scheme we first:  

1) Connect the given metric $g$ on $M^4$ with the standard metric of the direct product $\bar g$ on $\bar M^4=S^2_+\times S^2_-$ 
by some family of smooth metrics $g(\tau), 0\leq \tau\leq 1$; and then choose some partition $0<\tau_1<...<\tau_N<1$ such that the sequence of 
values ${\Cal H}(f^*, M^4, g(\tau_{k}))$ is sufficiently dense and starts with zero, see (5.1), (5.2) in the section~5 below. This is done with the sole purpose to prove that the minimal value ${\Cal H}(\Cal F)_{\alpha,\beta}$ is not achieved on the singular torus, i.e., for $(\alpha,\beta)$ from the boundary of $X$. If we would already know this - we would not need any sequence $g(\tau_k)$ and could start right with the given metric $g$.

2) Second: besides conditions (5.1) and (5.2) the sequence $g(\tau_k)$ should consists of the metrics satisfying some "transversality" conditions, which guarantee (as we wrote above) that under the given smooth family of corrections $f_\delta$ - zeroes of ${\Cal W}$ behave in a general way: for all but a finite number $\delta_i$ of instances of $\delta$ they are isolated and non-degenerated; while at each of the moment $\delta_i$ the pair of zeros of opposite indexes $\pm 1$ is born or annihilated. Such "bumpy metrics" could be constructed explicitly (in a standard way: taking partitions of unity and then adding small perturbations which make the holonomy - the operator of the "general" type). Or, this may be achieved by adding to the holonomy operator $I_{\gamma_{(\alpha,\beta,\phi,\psi)}}$ small "fake" holonomies $I^\epsilon_{(\alpha,\beta,\phi,\psi)}$ converging to zero as $\epsilon\to 0$.

\subhead 3.1. Corrections sequences  \endsubhead

To diminish the  functional we consider tori with positive contributions to the sum (1.3) and construct
"corrections" $\delta f$ around their extremals $\gamma^*$'s.

3.1. First step:  assume that for every extremal $\gamma$ we find a variational vector field $V$ as above in
2.2.1 in some neighborhood $U(\gamma)$ such that the derivative of ${\Cal H}(\gamma,f_\delta)$ on
$\delta$ for the variation of $f$ given by $V$ is negative. Using continuous dependence of ${\Cal H}$ and its
derivatives on $\gamma$ find possibly smaller neighborhood where such derivatives for other $\gamma'$ in this
smaller neighborhood still negative, say less then the half of the derivative for the extremal $\gamma$. Then
find a finite covering $\cup U_i$ of the extremal set by these neighborhoods and some partition of unity $\xi^i$
corresponding to this covering $U_i$ such that $\xi^i$ equals one in some smaller neighborhood of $\gamma\subset
U_i$. Next, consider the differential equation for $f_\delta$:
$$
{{d}\over{d\delta}} f_\delta = \xi^i V_i.
$$
From the tame  estimates above for the vector fields $V_i$ it follows the existence of the smooth solution of
the above equation on some small interval $0<\delta<\delta(f_0))$ - see [H]. As we already mentioned, our
corrections may be chosen to generate the Morse-type, or "shortening" deformations $T_\delta$ of the curves $T$
in the Grassmann bundle $G$ over $M^4$. In the same way as the Morse deformations gives in the limit a geodesic
line the above equation admits a solution $f_\delta$ converging to some $f^1$ with the extremals $\gamma(t)$
such that corresponding curves of planes $L(t)$ in $G_{\gamma(0)}$ are geodesics.

3.2. The second  step is to apply deformations of the types 2.2.2. Again, for every extremal $\gamma$ with
$L(t)$ being a geodesic define $V(\gamma)$ as in 2.2.2 above, find a small neighborhood $U(\gamma)$
of $\gamma$ where derivatives of the ${\Cal H}$ are sufficiently negative and proceed as above - continue $f^1$
over some small interval to a family of diffeomorphisms $f^1_\delta, 0\leq \delta \leq \delta^1$. Deformations
for the case when all $\gamma$'s on the given torus are extremals we described above in 2.2.4.

\medskip

These 3.1 and 3.2  are two building blocks of our Nash process. We do a sequence of successive corrections
applied to an initial diffeomorphism $f$ and then consecutive (intermediate) diffeomorphisms. The correction
process has an infinite sequence of "stages". "Each stage is like any other stage" and consists of two steps:
first we apply the (possibly) infinite number of steps 3.1 to obtain $f^1$ from $f$. Correspondingly, during the
first step in the number $(k+1)$ stage we construct $f^{k+1}$ from $f^k_{\delta^k}$. The second step is then
apply 3.2 to get $f^1_{\delta^1}$ from $f^1$, and correspondingly $f^k_{\delta^k}$ from $f^k$ during the second
step in the stage with number $k$.

Applying this double  induction we end up with the extremal map $f^*$ and some leaf ${\Cal F}^*$ stationary
under the above described variations for which it holds:
$$
{\Cal H}({\Cal F}^*)=0.  \qquad\qquad\qquad\qquad  {(\star)}
$$
Non-degeneracy of ${\Cal F}^*$ is a very important claim which we  establish next together with another crucial
fact: the index one of ${\Cal W}$ on $\partial X$ providing us with zeroes of index +1 inside $X$.

\medskip

\head 4. Non-collapsing. Boundary estimates \endhead

4.1. The crucial fact we  use below is that our functional ${\Cal H}(\Cal F_{\alpha,\beta})$ does not vanish on
the boundary of a parameter space: when some or both parameters $\alpha$ or $\beta$ converge by an absolute
value to $\pi /2$, i.e., when the tori-leaf $\Cal F_{\alpha,\beta}$ degenerates to the circle or a point - the
value of the functional ${\Cal H}(\Cal F_{\alpha,\beta})$ stays bounded from below by some positive constant
$K_{collapse}$\footnote{like $2\pi$} for arbitrary given diffeomorphism $f$. To see this we assume the contrary:
$$
{\Cal H}(\Cal F_{\alpha,\beta}) \to 0,  \tag 4.1
$$
and consider first the case  when only one of the parameters, say $\beta$ converge to $\pi/2$, while $\alpha$
stays bounded (by an absolute value) from it. Then $\Cal F_{\alpha,\beta}$ converges to the circle $f(\bar
l_{\alpha}(\phi))$ which we denote by $l(\phi)$. The second case is when $\Cal F_{\alpha,\beta}$ converge to the
point (one of the four fixed points of the $S^1$-action, singular points of $\Upsilon$).

\medskip

First case we divide further  into two subcases: first is when the operator $I_l$ is not identity. Then the
splitting of the tangent space to $M^4$ into the sum of invariant subspaces $\Pi$ and $N$ is correctly defined
for all curves close to $l$, e.g., from our family $\Upsilon$ on all tori in the considered family converging to
$l$. Hence, this splitting is defined over the union of $\Cal F_{\alpha,\beta}$ close to $l$, for instance, in
some tubular neighborhood of $l$. From (4.1) we see that all translations $I_\gamma$ for
$\gamma=\gamma_{\alpha,\beta}$ close enough to $l$ almost preserve tangent planes to $\Cal F_{\alpha,\beta}$,
which means that each tangent plane is close to corresponding $\Pi$ or $N$. Because the family of tangent planes
is continuous, while the distance from $\Pi$ to $N$ is a positive constant ($\pi/2$) this means that all tangent
planes of $\Cal F_{\alpha,\beta}$ are close, say, to $\Pi$, i.e., converge to the family of parallel planes
$\Pi(l(t))$ along $l$. Since $\dot l(t)$ belongs to $\Pi(l(t))$ (as the limit of vectors tangent to $\Cal
F_{\alpha,\beta}$), the plane $\Pi(l(t))$ is generated by $\dot l(t)$ and another unit vector field $Y(y)$
normal to $\dot l(t)$. Consider the "annulus" - the film $\Pi(s,t)=exp_{l(t)} sY(t)$. This film is normal to the
$N$-distribution along $l(t)$, therefore the normal projection on it of an arbitrary $\Cal F_{\alpha,\beta}$
sufficiently close to $l$ will be "an almost projection" along normal subspaces to $\Cal F_{\alpha,\beta}$,
i.e., a regular immersion of a tori into the "annulus" $\Pi(s,t)$ which is impossible.

Second subcase, when $I_l=id$,  is a little bit more involved. Now the splitting of the tangent spaces $T_pM^4$
into the sums of two invariant subspaces $\Pi$ and $N$ is no longer defined. Instead we may verify that (4.1) now
implies not only that the tangent spaces are almost parallel along $\gamma(t)$ from $\Upsilon$, but separately
each of the vector fields $X(t)$ and $Y(t)$ too. Indeed, the $\Cal F_{\alpha,\beta}$ is the boundary of the
image of the direct product under the smooth diffeomorphism $f$ of the product of the circle $\bar l_\alpha(t)$
and a small two-disk $\bar B^2(\epsilon)$ around some point (pole) in $S^2_-$. The tangent plane to this disk
goes under the differential of $f$ into the family of two-subspaces along $l(t)$ which we denote $N(t)$. If we
denote by $l_t(\psi)=f(\{\bar l(t)\}\times \partial \bar B^2(\epsilon))$ - "meridian" lines of $\Cal
F_{\alpha,\beta}$ then tangent vectors to these meridians are close to $N(t)$. Let $Y(t, \psi)$ be a unit
tangent vector to the meridian $l_t$ at $l_t(\psi)$, and $X(t,\psi)$ denotes, as before, the unit tangent vector
to the curve $\gamma\in\Upsilon$ from our one-dimensional foliation at the same point. Note, that by the
definition the curve $\gamma_\psi(t)$ issuing at $t=0$ from $l_0(\psi)$ satisfies $\gamma_\psi(t)=l_t(\psi+t)$;
so that (4.1) means that the tangent plane $T(t,\psi)$ generated by $X(t,\psi)$ and $Y(t,\psi)$ is almost
parallel along lines with coordinates $(t,\psi+t), 0\leq t\leq 2\pi$.

Combining the facts that 1)  under the parallel translation $I_\gamma(t)$ the plane $T(t,\psi+t)$ is almost
parallel due to (4.1); 2) this plane is generated by $X(t,\psi+t)$ and $Y(t,\psi+t)$; 3) that all $X(t,\psi+t)$
for different $\psi$ close to the same vector $\dot l(t)$, while the collection of $Y(t,\psi+t)$ for all $\psi$
tangent to meridians is close to the whole circle of unit vectors in the plane $N(t)$; 4) all parallel
transports $I_\gamma(t)$ along different $\gamma$ in $\Cal F_{\alpha,\beta}$ converge to the same operator
$I_l(t)$ - parallel transport along $l$; - we inevitably conclude that the parallel translation $I_\gamma(t)$
not only by (4.1) transport the plane $T(t,\psi+t)$ close to $T(0,\psi)$, but also that $I_\gamma(t) X(t,\psi+t)$
is close to $X(0,\psi)$ (which, in particular, means (in the limit as $\epsilon\to 0$) that $I_l(t)\dot
l(t)\equiv \dot l(0)$, or that $l(t)$ is a geodesic). Since the collection of all possible $Y(t,\psi+t)$ close
to $N(t)$, and $I_\gamma$, almost preserving $X(t,\psi+t)$ transport $T(t,\psi+t)$ close to $T(0,\psi)$
generated by $X(0,\psi)$ (close to $\dot l(0)$) and $Y(0,\psi)$ we conclude that in the limit $\epsilon\to 0$
all transports $I_\gamma$ preserve the family $N$: $I_\gamma N(t)=N(0)$; i.e., $N(t)$ is parallel along $l$.
Then the orthogonal complement $Z(t)$ to the three-space generated by $\dot l(t)$ and $N(t)$ gives a parallel
vector field along $l$.

Compare now $Y(0,\psi)$ with  the vector filed $Y^*(t,\psi)$ obtained by the parallel translations from
$l_t(\psi)$ along parallels of the tori. These "parallel's" are $l_\psi(t)=f(\bar l_\alpha(t)\times\{\bar
l_\beta(\psi)\})$. If $\tilde N(t,\psi)$ denotes the direction of the minimal geodesic connecting $l(t)$ with
$l_t(\psi)$ ("radius"), and $N(t,\psi)$ its projection on $N(t)$ then $N(t,\psi)$ is always transversal to the
projection of $Y(t,\psi)$ on the same $N(t)$. From which (and $I_l=id$ again) it is easy to conclude that for
each $0\leq t\leq 2\pi$ there exist a homotopy $H_t(\tau)$ from a normalized (unit) projections on $N(0)$ of the
vector field $Y^*(t,\psi)$ to $Y(0,\psi)$ continuously depending on $t$.

Now we get our contradiction as  follows: the unit projection of the vector field $Y(0,\psi)$ - tangent to the
meridian - on $N(0)$ gives a map of degree one of the meridian to the circle of unit vectors in $N(0)$. From the
above arguments it follows that it is homotopic to the map $\Delta  Y(\psi)$ obtained as follows: $\Delta
Y(\psi)$ is the parallel translation of $Y(0,0)$ first along $\gamma_0$ to $\gamma_0(\psi)$ and then ("back") to
$l_0(\psi)$ along "parallel" $l_\psi$. Note that by Ambrose-Singer theorem $\Delta Y(\psi)$ coincides with the
parallel transport $\tilde Y(\psi)$ of $Y(0,0)$ along the meridian to $l_0(\psi)$ up to the error proportional
to the area of the triangle with sides: $\gamma$, parallel $l_\psi$ and meridian, i.e., is small, or $\tilde
Y(\psi)$ is homotopic to $Y(0,\psi)$ for small $\epsilon$. Since, $\tilde Y(\psi)$ is null-homotopic for small
$\epsilon$ this contradicts the fact that $Y(0,\psi)$ has degree one.

\medskip

The second, last case left is  when the tori $\Cal F_{\alpha,\beta}$ converge to the point $P$. Let ${\Cal
F}(\epsilon)={\Cal F}_{\alpha(\epsilon),\beta(\epsilon)}$ be such a sequence, denote by $d(\epsilon)$ the
diameter of this tori ${\Cal F}(\epsilon)$, and re-scale the metric of $M^4$ in small neighborhood of $P$ by
$d^{-2}(\epsilon)$. In this new metric the family ${\Cal F}(\epsilon)$ has constant diameter 1 and converges to
the image of the direct product of two circles in flat $R^2\times R^2$ under the differential of $f$ at $\bar P,
f(\bar P)=P$. Then (4.1) would imply that ${\Cal H}(T^2)=0$ for the standard tori. We may see that this is not
true in yet another way. In the standard tori $T^2=\{(x,y,z,w)| x^2+y^2=1, z^2+w^2=1\}\subset R^4$ consider the
curve $\bar\gamma(t)$ issuing from $(0,1,0,1)$ in the direction $(1,0,1,0)/\sqrt{2}$, and find the curve
$\gamma(t)$ in ${\Cal F}(\epsilon)$ ($C^1$-)close to it. If ${\Cal H}(\gamma(t))$ is small, then the tangent
plane of ${\Cal F}(\epsilon)$ is almost parallel along $\gamma(t)$ implying that the derivative of the
$y$-coordinate of $\gamma$ is small, or $y(\gamma(t))$ stays close to 1. Therefore, $\gamma$ can not intersect
the "parallel" of ${\Cal F}(\epsilon)$ close to $(0,-1,cos(t), sin(t))$ in $T^2$. Which is impossible since the
homotopy class of $\gamma$ in $\pi_1({\Cal F}(\epsilon))$ is $(1,1)$ while the "parallel" is in $(0,1)$, i.e.,
they must intersect each other.

4.2. Now consider the vector  field ${\Cal W}$. For the standard direct product of two two spheres we easily see
that invariant planes $\Pi_{\alpha,\beta}$ and $N_{\alpha,\beta}$ are generated by the tangent planes to the
factors of this direct product; and that ${\Cal W}({\alpha,\beta})=(sin(\alpha),sin(\beta))$\footnote{again, by Gauss-Bonnet theorem.}.
As we've just seen, under the change of the metric the maximum of ${\Cal H}(\gamma_\psi)$ on $\psi$ does not go to zero, which
implies that some ${\Cal W}(\psi)$ also stay bounded away from the origin. Since all parallel transports
$I_{\gamma_\psi}$ for collapsing tori (thin or converging to a point) converge to each other this implies that
the integral ${\Cal W}$ of ${\Cal W}(\psi)$ also is bounded from zero, or that the vector field ${\Cal W}$ does
not vanish on the boundary. Hence, it is homotopic to the initial one, or has index one.
$$
ind({\Cal W})_{| \partial X} =1. \tag 4.2
$$

\medskip

\head 5. The Hopf conjecture \endhead

Connect the given metric $g$ on $M^4$ with the standard metric of the direct product $\bar g$ on $\bar
M^4=S^2_+\times S^2_-$ by some family of smooth metrics $g(\tau), 0\leq \tau\leq 1$. For $\tau=0$ the metric
$g(0)=\bar g$ admits the "silent" $\Cal F^*$ leaf with ${\Cal H}(id)=0$ - it is the product of equators of the
factors $S^2_+$ and $S^2_-$ - the totally geodesic torus. As we already mentioned when defining ${\Cal H}$, for
the general family of metrics ${\Cal H}(f^*, M^4, g(\tau))$ continuously depends on $\tau$. Hence, for some
small $\tau_1$ and all $0<\tau\leq \tau_1$ it holds ${\Cal H}(f^*, M^4, g(\tau))\leq K_{collapse}$, where
$K_{collapse}$ is the constant from the previous chapter, which implies that corresponding $\Cal F^*$ is a
regular tori if it exists. Further, find a partition $0<\tau_1<...<\tau_N<1$ such that
$$
|{\Cal H}(f^*, M^4, g(\tau_{k+1}))- {\Cal H}(f^*, M^4, g(\tau_{k}))|\leq K_{collapse}/2 \tag 5.1
$$
and
$$
|{\Cal H}(f^*, M^4, g(1))- {\Cal H}(f^*, M^4, g(\tau_{N}))|\leq \epsilon \to 0. \tag 5.2
$$
and such that $g(\tau_k), k=1,...,N$ satisfy the standard assumptions. Applying the Nash process we can establish
by induction that ${\Cal H}(f^*, M^4, g(\tau_k))=0$ for all $k$, and, taking $\epsilon\to 0$ that  ${\Cal H}(f^*, M^4, g(1))=0$ too.
For each $\epsilon$ we have the regular tori $\Cal F^*$ of some flag foliation $(\Upsilon, \Cal F)$ on $(M^4,g(\tau_N))$
under some diffeomorphism $f^*$ such that ${\Cal H}(\Cal F^*)=0$. By Gauss theorem its extrinsic curvature vanish, i.e.,
the integral of the inner (Gaussian) curvature of ${\Cal F}^*$ equals the integral of the sectional curvature of
$M^4$ tangent to the tori. Since, by the Gauss-Bonnet theorem the integral of its inner curvature of any tori equals
zero, the sectional curvature of $(M^4,g(\tau_N))$ vanishes somewhere on ${\Cal F}^*$. Hence, by continuity the same is true for
$(M^4,g)$, and the Hopf hypothesi follows.

\medskip


\head Appendix~1. Flag Distributions Functionals \endhead

A1.1. The text above is an  application of a more general concept - flag distributions functionals. Originally
they were suggested for the construction introducing the metric and the Hodge-star operator on any manifold
endowed with two flag distributions. There the non-integrability of distributions plays an important role. For
our particular purposes though we studied here a manifold endowed with only one flag distribution which is
integrable - a flag foliation; i.e. we used "half" of the construction, but already "charged" with some metric.
The flag foliation on a Riemannian manifold $(M^n,g)$ is defined when at every point $x$ we have a flag of
vector spaces $V_1 < V_2 < ... < V_k$ in the tangent space $T_xM$ such that each distribution $V_i$ is the
tangent one to some foliation $F_i$ of the manifold $M$. We say that the collection $\{F_i, i=1,...,k\}$ is the
$k$-flag of foliations on $M$. In general we allow foliations to be singular; hence, the dimension of a
particular distribution $V_i$ is not constant, but only semi-continuous. To give the idea of a flag foliation
functional we restrict ourself to the simplest case of a flag of two foliations as above: the pair $(\Upsilon,
\Cal F)$ of a one and two-dimensional foliations on $M$.

In this case the constructions of some (other than ${\Cal H}$) functionals  may proceed as follows.

\medskip

\subhead 1. Relative, integral  and module $\Cal H$-functionals \endsubhead First, for an arbitrary curve
$\gamma(s)$ of the family $\Upsilon$ denote by $\bar N_\gamma$ its normal bundle in $M$, by $N_\gamma$ the
restriction of the normal bundle of $\Cal F$ to $\gamma$, and by $G_\gamma$ - the Grassmann bundle of $(codim
\Cal F)$-subspaces in $\bar N_\gamma$. By $\tilde\gamma(s)$ we denote the natural lift of $\gamma(s)$ into
$G_\gamma$ where $\tilde\gamma(s)$ equals $N_\gamma(s)$ - is the normal to $\Cal F$ subspace at the point
$\gamma(s)$. Denote by $\nabla$ the Riemannian connection of $M$, and by $\nabla^*_\gamma$ (or simply by
$\nabla^*$ when it does not lead to a confusion) the component of $\nabla$ normal to $\gamma$. In other words:
for an arbitrary vector field $Y(s)$  along $\gamma$
$$
\nabla^* Y(s)=(\nabla_{\dot\gamma} Y(s))^\perp,  \tag A1.1
$$
where $(...)^\perp$ denotes the  component of the vector normal to $\dot\gamma$. This connection defines a
metric in $G_\gamma$ as follows: for an arbitrary curve $\tilde L(s)$ in $G_\gamma$ (i.e., the curve of
$codim(\Cal F)$-subspaces in the curve $N_\gamma$ of $codim(\gamma)$-subspaces) connection $\nabla^*$ induces
connection on $G_\gamma$, which we denote also by $\nabla^*$ and such that $\nabla^*\tilde L(s)$ is tangent to
$G_{\gamma(s)}$. The norm of this vector in the natural metric on the Grassmannian $G_{\gamma(s)}$, say, coming
from the Riemannian submersion from $O(n)$, equals the norm of $\dot{\tilde L(s)}$. In the case we consider in
this note -- $\gamma$ is a curve in some leaf $F$ of the two-dimensional foliation $\Cal F$ of the
four-dimensional $M$. Then $G_\gamma$ is a bundle over $\gamma$ of two-subspaces in the bundle $N_\gamma$ over
$\gamma$ of three-subspaces; i.e., subspaces normal to $\gamma$ in $M$. Let  $Y(s)$ be a unit vector field
normal to $\gamma$ and tangent to $F$, and $W(s),Z(s)$ be some unit and normal to each other vector fields along
$\gamma$ normal to $F$. Then $\tilde\gamma(s)$ is the curve of two-subspaces generated by $\{W(s),Z(s)\}$ in the
three-subspace $N_{\gamma(s)}$ generated by $\{Y(s),W(s),Z(s)\}$. Since $\tilde\gamma(s)$ is uniquely determined
by the vector $Y(s)$ in our case the norm of $\dot{\tilde L(s)}$ simply equals the norm of
$\nabla_{\dot\gamma(s)} Y(s)$:
$$
\|\dot{\tilde L(s)}\|^2=\|\nabla^*_{\dot\gamma(s)} Y(s)\|^2=(\nabla_{\dot\gamma(s)}Y(s),W(s))^2+(\nabla_{\dot\gamma(s)}Y(s),Z(s))^2. \tag A1.2
$$
Take for a particular curve $\gamma(s), 0\leq s\leq 1$ from $\Upsilon$
$$
t^{MA}_{VA}(\gamma, F)=  \int\limits_0^1 \|\dot{\tilde L(s)}\|^2 ds,  \tag A1.3
$$
and
$$
st^{MA}_{VA}(M,(\Upsilon,{\Cal F}))= sup\{ t^{MA}_{VA}(\gamma) | \gamma\in \Upsilon\}. \tag A1.4
$$
where $t^{MA}_{VA}(\gamma)$ is $t^{MA}_{VA}(\gamma, F)$ for the surface $F$ from $\Cal F$ which contains $\gamma$.

We say that the flag foliation  $(\Upsilon,{\Cal F})$ admits an invariant transversal measure if the foliation
$\Cal F$ admits an invariant transversal measure $d\mu_{\Cal F}$ and almost $d\mu_{\Cal F}$-everywhere the
sub-family $\Upsilon_F$ of our curves $\gamma$ in the given leave $F$ of $\Cal F$ is also measurable. Then the
integral over all $\Upsilon$ is correctly defined. In this case we take
$$
t^{MA}_{VA}(M,(\Upsilon,\Cal F))= \int\limits_{\Cal F} \int\limits_{\Upsilon_F} t^{MA}_{VA}(\gamma,F) d\Upsilon_F d\mu_{\Cal F}. \tag A1.5
$$
{Let ${\| {\gamma} \|}_{L^2}$ be the energy of $\gamma$. If to  denote $t^{MA}_{VA}(\gamma,
F)=\|\tilde\gamma\|^2_{L^2}$ then instead of (A1.3) another obvious choice is: $et^{MA}_{VA}(\gamma,
F)=\|\gamma\|^2_{L^2} + \|\tilde\gamma\|^2_{L^2}$. Accordingly, we obtain another functionals
$s't^{MA}_{VA}(M,(\Upsilon,\Cal F))$ by taking instead or sup-norm as in (A1.4) over all $\gamma$; or
$e't^{MA}_{VA}(M,(\Upsilon,\Cal F))$ - an integral ($L^2$)-norm as in (A1.5).

Other interesting possibility might be the "conformal" functional
$$
ct^{MA}_{VA}(M,(\Upsilon,\Cal F)) = sup \{ \Lambda(\Upsilon,\rho)| \rho \in L^2(F), \|\rho\|_{L^2}=1 \}, \tag A1.6
$$
where
$$
\Lambda(\Upsilon,\rho) = inf \{ \int \|\dot{\tilde\gamma(s)} \rho(\gamma(s)) ds  \| | \gamma\in\Upsilon\}.
$$}

\medskip

A1.2. Even more interesting is  the "non-commutative"  $nt^{MA}_{VA}(M,(\Upsilon,\Cal F))$ one (a'la Connes),
when $\Upsilon$ is understood as an operator (i.e., an element of the von Neumann algebra $W^*(M,\Cal F)$ of
bounded random operators acting (as an integral above) on sections of the Hilbert bundle $H_x=L^2(\Cal F(x)),
x\in {\Cal F}(x)$, see [C].\footnote{ {\bf Foliation immersions.} Note that for more general than
diffeomorphisms maps $f: \bar M\to M$ the definition above may be still possible: Let $(\bar M,(\bar\Upsilon,
\bar {\Cal F}))$ be a manifold with a 2-flag foliation, and $f:\bar M\to M$ some differentiable map. Then the
image of the flag foliation is not necessary a flag foliation due to possible degenerations of $f$. But if $f$
is regular enough, say its restriction on every leaf $\bar F$ of $\bar {\Cal F}$ is an immersion $f_{|\bar F}:
\bar F\to F$ then $t^{MA}_{VA}(M,\gamma,F)$ is defined for every curve $\gamma=f(\bar\gamma)$ for
$\bar\gamma\in\bar\Upsilon$, and correspondingly -- the supremum of all $t^{MA}_{VA}(M,\gamma,F)$ as in (4').
The obtained number we denote by $t^{MA}_{VA}(M,f)$.}

Easy to see that our ${\Cal W}$  is a remnant ("vector trace") of some operator from $W^*(M,\Cal F)$ of this type.
When some $M$ is endowed with a flag foliation $(\Upsilon,\Cal F)$ such that the moduli space $X=M / {\Cal F}$ is
not Hausdorff (say, considering the "rank problem" which is the kind of the non-commutative Hopf conjecture) this
"non-commutative" approach seems very naturale.

\medskip

\subhead 2. Further applications. Dutch category \endsubhead

A2.1. The text above may be  understood as follows: a one-dimensional foliation $\Cal \Upsilon$ defines the
family of parallel transports $I_\gamma, \gamma\in\Upsilon$ which in turn defines the distribution of two
orthogonal plane fields $\Pi$ and $N$ - invariant subspaces of $I_\gamma$. It defines an almost complex
structure $J_\Upsilon$ on $M^4$ - rotations on $\pi/2$ about $\Pi$ and $N$. What we proved, actually, was that
there exists a flag-foliated diffeomorphism $f^*: (\bar M,\bar\Upsilon,\bar\Cal F)\to (M,\Upsilon,\Cal F)$
admitting the ("silent") leaf $\Cal F^*$ which is

1) "$\Upsilon$-holomorphic"\footnote{another  connection to the complex geometry may be through the identity
between the Grassmannian of two-planes in $R^4$ and lines in $CP^3$, i.e., $so(4)\equiv su(2)$}, i.e., tangent
to the one of the distributions of invariant subspaces of the almost complex structure $J_\Upsilon$, and

2) the almost complex  structure is "$\Upsilon$-parallel" along this leaf: $\nabla_{\dot\gamma} J_\Upsilon
\equiv 0$.

3) non-collapsing - in  some important cases we can be sure that $\Cal F^*$ is not degenerate.

Then our functional  $\Cal He$, appears, count these "$\Upsilon$-holomorphic" leaves with the help of the Euler
co-cycle $\Xi=\sum ind(x)$ acting on the fundamental class of the moduli space $X^2=M^4/T^2$ with weight ${\Cal
H}$. The whole text above then provides us with the existence of the "$\Upsilon$-holomorphic" leaf. By
definition it is defined on the Teichm{\" u}ller space ${\Cal T}(M)={\text Metrics}(M^4)/{\text Diff}(M^4)$. We
may speculate a little bit more on how further "complexification" goes: if we think of $W(\alpha,\beta)$ as a
new complex coordinate (instead of $\{\alpha,\beta\}$) in the moduli space $X$, and assume that for some
diffeomorphism $f$ the corresponding ${\Cal H}({\Cal F}_{\alpha,\beta})$ is holomorphic in this new
coordinates\footnote{should we be looking for such pairs $({\Cal W,H})$, actually?} then our functional is given
by
$$
{\Cal He}(M,f,g) = {{1}\over{2\pi}}\int\limits_{\partial X} {{\Cal H(W) dW}\over{\Cal W}}. \tag A2.1
$$

A2.2. Possibly (for "bumpy" metrics) we may add many more objects of this type by changing the homotopy class of
the wirings $\Upsilon$ of the foliation $\Cal F$ or taking another torus action on $M^4$, and then study how this
collection - $(\Upsilon,{\Cal F})$-"etazherka"\footnote{Wh\'{a}t-n\`{o}t} - behaves under the change of the class of the
metric in ${\Cal T}(M)$ (which will gives as "the functor" from the category of paths in ${\Cal T}(M)$ to our
own). This structures may play the role of the cut-sets to double in order to decompose the manifold into simpler pieces.

Note also that there exists an  obvious duality between the classes of flag foliations and open book decompositions: in
the first one the disjoint union of lower dimensional leaves is the leaf of a higher dimension, while in the
second the intersection of disjoint leaves of a higher dimension is the leaf of a lower dimension. This duality
may be extended further to corresponding form-descriptions of flag foliations and book decompositions. If we
call flag foliations "closed books" then $(\Upsilon,{\Cal F})$-"etazherka" would be an appropriate place to
store them. It is also possible to composite both flag foliations and open book decompositions together, say considering
the open book whose pages are flag foliations; etc.



\enddocument
\bye